\documentclass[microtype]{gtpart}
\usepackage{indentfirst,epic,url,graphics,needspace}
\usepackage{hyperref}

\title{Simplicial complexes with lattice structures}

\author{George M. Bergman}
\givenname{George}
\surname{Bergman}
\address{Department of Mathematics\\
University of California\\
Berkeley, CA 94720-3840\\USA}
\email{gbergman@math.berkeley.edu}
\urladdr{http://math.berkeley.edu/~gbergman}

\arxivreference{1602.00034}
\arxivpassword{pf42m}

\newtheorem{theorem}{Theorem}
\newtheorem{lemma}[theorem]{Lemma}
\newtheorem{corollary}[theorem]{Corollary}
\newtheorem{proposition}[theorem]{Proposition}
\newtheorem{definition}[theorem]{Definition}
\newtheorem{question}[theorem]{Question}
\newtheorem{convention}[theorem]{Convention}

\renewcommand{\r}{\mathrm}
\newcommand{\dotline}[2]{\put(#1){\circle*{2}}\qbezier(#1)(#1)(#2)}
\newcommand{\ho}{(0,1]} % half-open interval
\mathchardef\hy="2D

\DeclareRobustCommand{\smallcoprod}{\begin{picture}(10,6)
\thicklines
\put(1.8,.6){\line(1,0){6.6}}
\put(3.6,.6){\line(0,1){6}}
\put(6.6,.6){\line(0,1){6}}
\end{picture}}

\begin{document}

\begin{abstract}
If $L$ is a finite lattice, we show that there is a natural
topological lattice structure on the geometric realization
of its order complex $\Delta(L)$ (definition recalled below).
Lattice-theoretically, the resulting object
is a subdirect product of copies of $L.$
We note properties of this construction and of some
variants, and pose several questions.
For $M_3$ the $\!5\!$--element nondistributive
modular lattice, $\Delta(M_3)$ is modular, but
its underlying topological space does not admit a structure
of distributive lattice, answering a question of Walter Taylor.

We also describe a construction of ``stitching together'' a
family of lattices along a common chain, and note
how $\Delta(M_3)$ can be regarded as an example of this construction.
\end{abstract}
\maketitle
% - - - - - - - - - - - - - - - - - - - - - - - - - - - - - -

\section{A lattice structure on \texorpdfstring{$\Delta(L)$}{Delta(L)}}\label{S.constr}

I came upon the construction studied here
from a direction unrelated to the
concept of order complex; so I will first motivate
it in roughly the way I discovered
it, then recall the order complex construction, which turns
out to describe the topological structures of these lattices.

\subsection{The construction}\label{S2.constr}
The motivation for this work comes from Walter Taylor's paper \cite{WT},
which examines questions of which topological spaces --
in particular, which
finite-dimensional simplicial complexes -- admit various sorts
of algebraic structure, including structures of lattice.
An earlier version of that paper asked
whether there exist spaces which admit structures
of lattice, but not of distributive lattice.
(More precisely, it asked whether there are spaces admitting
structures of modular lattice but not of distributive lattice,
and whether there are spaces admitting
structures of lattice but not of modular lattice.
But the above simplified version is enough to
motivate what follows.)

The most familiar examples of nondistributive lattices are
finite lattices such as $M_3=\!$
\raisebox{2.5pt}[12pt][7pt]{ % M_3
\begin{picture}(15,20)
\dotline{7,-8}{7,0}
\dotline{7,0}{7,8}
\dotline{7,8}{0,0}
\dotline{0,0}{7,-8}
\dotline{7,8}{14,0}
\dotline{14,0}{7,-8}
\end{picture}}\,
and $N_5=$%
\raisebox{2.5pt}[12pt][7pt]{ % N_5
\begin{picture}(15,20)
\dotline{6,-8}{0,-3.5}
\dotline{0,-3.5}{0,3.5}
\dotline{0,3.5}{6,8}
\dotline{6,8}{12,0}
\dotline{12,0}{6,-8}
\end{picture}},
whose underlying sets, looked at as discrete topological
spaces, certainly also admit structures of distributive lattice
(e.g., they can be rearranged into chains).
But we may ask whether, starting with such a finite lattice $L,$
there is some way of building from it
a more geometric sort of lattice, whose
underlying topological space has a distinctive geometry
which perhaps precludes a distributive topological lattice structure.

As a first attempt, one might identify each $x\in L$ with
the $\![0,1]\!$--valued function on $L$ having value $1$
at $x$ and $0$ elsewhere, and try to extend the lattice structure
on $L$ to convex linear combinations of these functions.
However, there is no evident way of defining meets and joins of such
linear combinations so that they extend the given operations on
$L,$ while continuing to satisfy the lattice identities.

But suppose, instead, that we identify each element $x$ of our
finite lattice $L$ with
the $\![0,1]\!$--valued function having value $1$ on the principal
ideal ${\downarrow}\,x=\{y\in L\mid y\leq x\}$ generated by $x,$
and $0$ elsewhere.
I claim that these functions
belong to a family which does have a natural lattice structure:
the space $F(L)$ of all $\![0,1]\!$--valued functions $f$ on $L$
such that for each $t\in [0,1],$ the set
\begin{equation}\begin{minipage}[c]{25pc}\label{d.f_t}
$f_t\ =\ \{x\in L\mid f(x)\geq t\}$
\end{minipage}\end{equation}
is a principal ideal of $L;$ and that the characteristic functions
we started with form a sublattice of $F(L)$ isomorphic to $L.$

Indeed, let us order functions $f$ with the above property by pointwise
comparison, writing $f\leq g$ if $f(x)\leq g(x)$ for all $x\in L.$
A greatest lower bound of any $f,\,g\in F(L)$
is given by their pointwise infimum $f\wedge g,$
since for each $t\in[0,1],$ the set $(\,f\wedge g)_t$ is
$f_t\cap g_t,$ an intersection of principal
ideals, and hence again a principal ideal of~$L.$

That $f$ and $g$ have a least upper bound can be seen in various ways.
On the one hand, one can define $f\vee g$ to be the pointwise infimum
of all upper bounds of $f$ and $g$ in $F(L),$ and note that for
all $t\in[0,1]$ this makes $(\,f\vee g)_t$ the intersection of a set of
principal ideals, which, since $L$ is finite, is again a
principal ideal.
Hence $f\vee g\in F(L),$ and by construction it
will be the least member of $F(L)$ majorizing both $f$ and $g.$
Alternatively, one can guess that the desired function
$f\vee g$ on $L$ will have the property that for each $t\in[0,1],$ the
set $(\,f\vee g)_t$ is the join of $f_t$ and $g_t$ as ideals of $L,$
and verify that this condition indeed
determines an element $f\vee g\in F(L),$ and
gives the desired least upper bound.
Finally one can, for each $x\in L,$ take $(\,f\vee g)(x)$
to be the largest $t\in[0,1]$ such that $x$ is majorized by the
join of an element $y$ such that $f(y)\geq t$
and an element $z$ such that $g(z)\geq t,$ and again verify that
the function $f\vee g$ so defined has the desired properties.
Since least upper bounds are unique when they
exist, all these constructions give the same operation on $F(L).$

So $F(L),$ regarded as a partially ordered set under coordinatewise
comparison, forms a lattice, and we
see that the elements of $F(L)$ that are $\!\{0,1\}\!$--valued
form a sublattice isomorphic to $L.$
In \S\ref{S2.M3} we shall see that when applied to $M_3,$
this construction gives a topological modular lattice whose
underlying topological space does not admit a distributive lattice
structure, and so indeed answers one part of the question mentioned.

\subsection{The order-complex construction}\label{S2.Delta}
For $F(L)$ defined as above, any $f\in F(L)$ can be written as a
convex linear combination of the characteristic functions of a
{\em chain} of principal ideals, namely, the sets $f_t$ $(t\in[0,1]).$
The characteristic functions in such a chain are linearly
independent, and the principal ideals in terms of which $f$ is
so expressed are unique if we require that each characteristic
function appear with nonzero coefficient.
For each such chain of $n$ characteristic functions, the
$\!n\!$--tuples of coefficients that can be applied to them
to get a convex linear combination (with coefficients now allowed
to be zero) are the $\!n\!$--tuples of real numbers in $[0,1]$ which sum
to $1,$ and these form, geometrically, an $\!(n{-}1)\!$--simplex.
The observation about unique representations with nonzero coefficients
shows that though a given point of $F(L)$ may lie in the
simplices corresponding to more than one chain, it will
lie in the interior of only one.
Thus, $F(L)$ has the form of a simplicial complex,
with simplices corresponding to the chains of principal ideals of $L;$
equivalently, to the chains of elements of $L.$

The geometric construction just described is an instance
of a standard concept, the {\em order complex} of a poset $P.$
This is the simplicial complex $\Delta(P)$ having an $\!n\!$--simplex
for each chain of $n+1$ elements in $P,$ and such that the faces of
the simplex determined by a chain $C$ are the simplices determined
by the subchains of $C.$

The term ``simplicial complex'' is commonly applied both to
(i)~an abstract object given by a set $V$ (``vertices'') and
a family $D$ of subsets of $V$ (``simplices''; an $\!(n{+}1)\!$--element
set being called an $\!n\!$--simplex), such that $D$ is
closed under passing to subsets (faces of simplices), and also to
(ii)~the {\em geometric realization} of such
an abstract simplicial complex: a topological space
in which the abstract simplices are replaced by
subspaces -- a point for each $\!0\!$--simplex, a line segment for each
$\!1\!$--simplex, a triangle for each $\!2\!$--simplex, and so forth
-- with the appropriate face-relations among them.
The simplicial complex, in this sense, is the union of these simplices.
So what we have shown is that the {\em geometric realization} of the
order complex of a lattice has a natural lattice structure.

In the literature, the geometric realization of an abstract simplicial
complex $D$ may either be denoted by a
symbol such as $|D|$ or $||D||,$ or by the same symbol $D,$ allowing
context to determine which is meant (cf.\ \cite{MW}, sentence
beginning at bottom of p.\,5).
Here we shall always understand simplicial complexes in
the geometric sense, unless they are specified as ``abstract''.
We will occasionally say ``geometric realization'' for emphasis.

From now on, in place of the notation $F(L)$ used above,
I will use the standard notation, $\Delta(L).$

The construction of the order complex of a partially
ordered set $P$ appears to have been introduced in 1937 by
P.\,Alexandroff \cite[\S2.1]{PA}, who called it the ``barycentric
subdivision'' of $P,$ because the barycentric subdivision, in the
classical sense, of
a simplicial complex $D$ is given by the order complex of the
partially ordered set of simplices of $D,$ ordered by inclusion.
(In \cite[p.\,7]{MW} the operation of barycentric subdivision
on simplicial complexes is described, more precisely,
as the composite of the construction taking a simplicial
complex to the partially ordered
set of its simplices, and the construction
taking a partially ordered set to its order complex.)
The concept seems to have been revived, or perhaps rediscovered,
by Jon Folkman \cite[next-to-last sentence before Theorem 3.1]{JF},
and is now a standard tool \cite{MW}, \cite{JW}.
What may be new in the approach
of this note is the representation of
the points of these complexes as convex linear combinations of
characteristic functions of principal ideals.
As noted, the set of functions that can be so described has a natural
partial ordering by pointwise comparison, and when $P$ is a lattice,
the resulting poset has the lattice structure discussed above, which
would otherwise be hard to see.

For $L$ a finite lattice, let us define for
each $t\in[0,1]$ the function $h_t\,{:}\ \Delta(L)\to L,$
associating to each $f\in\Delta(L)$ the generator $h_t(\,f)$ of the
ideal $f_t=\{y\in L\mid f(y)\geq t\}.$
It follows from our description of the lattice operations
on $\Delta(L)$ that each $h_t$ is a lattice homomorphism;
and clearly the family of these functions as $t$ ranges
over $[0,1]$ separates points of $\Delta(L).$
\vspace{.2em}

In the next theorem we summarize the above observations, putting them
initially in the context of a general finite poset $P.$
We will not distinguish notationally between a poset or lattice
and its underlying set; so $[0,1]^P$ will mean the set
of all $\![0,1]\!$--valued functions on the set $P.$
We will write ``principal ideal'' for ``principal downset''
in speaking about posets, so as to use the same language
for posets and lattices.
The reason why the condition that $f_t$ be a principal ideal
of $P$ is imposed, below, only for $t\in\ho,$
rather than for all $t\in[0,1]$ as above,
will be noted in \S\ref{S2.01}.

\begin{theorem}\label{T.main}
If $P$ is a finite poset, then the geometric realization of
its order complex, $\Delta(P),$ can be identified with the set of all
functions $f\in[0,1]^P$ having the property that for each
$t\in\ho,$ the subset $f_t=\{x\in P\mid f(x)\geq t\}$ is
a principal ideal of $P.$
For such $f$ and $t,$ we shall write $h_t(\,f)$ for the element
of $P$ such that $f_t={\downarrow}\,h_t(\,f).$
Regarding $\Delta(P),$ so defined, as partially ordered by pointwise
comparison, each map $h_t$ is an isotone map $\Delta(P)\to P,$
and this family of maps separates points of $\Delta(P).$

If the finite poset $P$ has pairwise meets
\textup{(}greatest lower bounds\textup{)} of elements,
then so does $\Delta(P).$
These are given by pointwise infima of functions $P\to [0,1];$
thus, they may be described by the formula
\begin{equation}\begin{minipage}[c]{25pc}\label{d.wedge}
$(\,f\wedge g)_t=f_t\wedge g_t=f_t\cap g_t,\quad(\,f,g\in\Delta(P),$
$t\in\ho).$
\end{minipage}\end{equation}
This operation $\wedge$ is continuous on $\Delta(P),$
and is respected by the maps $h_t.$

Similarly, if $P$ has pairwise joins
\textup{(}least upper bounds\textup{)} of elements,
then so does $\Delta(P).$
These may be described by the formula
\begin{equation}\begin{minipage}[c]{25pc}\label{d.vee}
$(\,f\vee g)_t=f_t\vee g_t,\quad(\,f,g\in\Delta(P),$ $t\in\ho).$
\end{minipage}\end{equation}
where $f_t\vee g_t$ denotes the join as principal ideals.
The resulting operation $\vee$ on $\Delta(P)$ is,
again, continuous and respected by the maps $h_t.$
If in addition to pairwise joins, $P$ has a join of the empty family,
i.e., a least element $0,$ then the elements of $\Delta(P)$ can be
described as those maps $f\,{:}\ P\to[0,1]$ which satisfy the identities
\begin{equation}\begin{minipage}[c]{25pc}\label{d.vee=wedge}
$f(0)=1$\quad and\quad $f(x\vee y)=f(x)\wedge f(y)$ $(x,y\in P).$
\end{minipage}\end{equation}

Thus, if $L$ is a finite lattice, then $\Delta(L),$ ordered
as above by pointwise comparison, is a
topological lattice, which as an abstract lattice is a subdirect
product of copies of $L.$
Again, its elements can be described as
all maps $f{:}\ L\to[0,1]$ satisfying~\eqref{d.vee=wedge}.
\end{theorem}

\begin{proof}[Sketch of proof]
Everything has been covered in the preceding remarks except
the sentence containing~\eqref{d.vee=wedge} (and its echo in the
final sentence).
To see that $f(x\vee y)=f(x)\wedge f(y)$ holds when
$f\in\Delta(P)$ and $x,y\in P,$
let $s=f(x\vee y)$ and $t=f(x)\wedge f(y).$
Then the facts that $f_s$ and $f_t$ are principal downsets give
the relations $\leq$ and $\geq$ respectively.
The identity $f(0)=1,$ and the fact that these two identities
imply membership in $\Delta(P),$ are straightforward.
\end{proof}

(One may ask: if $P$ is a finite partially ordered set merely having
pairwise meets, why can't we get joins $f\vee g$ in $\Delta(P)$ by
one of the constructions described earlier,
namely, as the pointwise infimum of
all upper bounds of $f$ and $g$ in $\Delta(P)$?
The answer is that if we don't assume $P$ also has joins,
the set of upper bounds of $f$ and $g$ in $\Delta(P)$ may be empty.)

Note that the maps $h_t$ in the above theorem,
since they send the connected space $\Delta(P)$ to the
discrete space $P,$ are necessarily discontinuous; hence the
resulting map $\Delta(P)\to P^{\ho}$
defining our subdirect product structure is also discontinuous
with respect to the natural topology on $P^{\ho}.$
More on this in \S\ref{S2.subdirect}.

\subsection{Notes on \texorpdfstring{$0$ and $1\!$}{0 and 1}}\label{S2.01}
In the construction of $\Delta(P)$ in Theorem~\ref{T.main}, we
required $f_t$ to be a principal ideal for $t\in\ho,$
but not for $t=0.$
This is because every $\![0,1]\!$--valued function $f$
on $P$ satisfies $f_0=P;$ hence if $P$ has no greatest element,
the set of functions $f{:}\ P\to[0,1]$ such that
$f_0$ is a principal ideal is empty; so $\Delta(P)$ would be
empty if we required its elements to have that property.
On the other hand, for posets $P$
{\em with} greatest element, such as finite lattices,
the condition that $f_0=P$ be a principal ideal is vacuous,
so its omission makes no difference.
If we were only considering lattices, it would be natural to
word our condition as saying
$f_t$ is a principal ideal for all $t\in[0,1];$
but since we will be proving many of our results for general
posets, we require this only for $t>0.$

At the other end, notice that for $f\in\Delta(P),$
since $f_1$ is a principal ideal, it is nonempty.
Thus, as $P$ is finite, $f$ must have the value $1$
at some minimal element of $P.$
If $P$ has a least element $0_P$ -- in particular, if it is lattice --
this tells us that $f(0_P)$ is automatically $1.$
So in discussing the values of a function $f\in\Delta(L)$ for
$L$ a lattice, we may take for granted the condition
$f(0_L)=1,$ and discuss $f$ in terms of its values on $L-\{0_L\}.$

\subsection{Some conventions}\label{S2.conventions}
Let us make precise some terms we have already used.

\begin{convention}\label{C.nonempty}
In this note, all lattices, semilattices and partially
ordered sets will be assumed nonempty.
\end{convention}

This is mainly so that finite lattices are guaranteed to have
a least element $0.$
(Alternatively, we could have supplemented the lattice operations with a
zeroary join operation, giving $0$ as the
join of the empty family, alongside pairwise joins.
Still another possibility would have been to omit
Convention~\ref{C.nonempty},
and allow the empty lattice along with the others,
noting, for instance, that $\Delta(\emptyset)=\emptyset,$
since the unique $\![0,1]\!$--valued function $f$ on the
empty set does {\em not} have the property that the
sets $f_s$ are principal ideals -- since the empty lattice has no
principal ideals.
But I chose to avoid the distractions that
that special case would have entailed.)

\begin{convention}\label{C.topological}
In this note, a {\em topological lattice} means
a Hausdorff topological space with a lattice structure under which the
lattice operations are continuous.
\end{convention}

We remark that some natural-seeming topologies on
lattices fail to satisfy these conditions.
For instance, the lattice of vector subspaces of $\mathbb{R}^2$
may be regarded as composed
of a circle, comprising the $\!1\!$--dimensional subspaces,
and two isolated points, the trivial subspace and the whole space.
But the meet and join operations are discontinuous on the circle,
since their outputs jump from the isolated points
back to the circle whenever their two arguments fall together.
We also note that in the
literature, a partially ordered set is sometimes given
the non-Hausdorff ``Alexandrov (or Alexandroff) topology'' introduced
in \cite{PA}, in which the open subsets are the downsets.
But we will not consider such topological spaces here.

It will be convenient to make a notational identification:

\begin{convention}\label{C.downarrow}
If $x$ is an element of a partially ordered set $P,$ then
the symbol ${\downarrow}\,x,$ denoting the principal ideal
of $P$ generated by $x,$ will also be used, in discussing
$\Delta(P),$ for the characteristic function of that ideal.
\end{convention}

One more remark on usage in the literature: an abstract
simplicial complex is in some works, such as \cite{MW},
assumed closed under taking {\em all} subsets of its members,
while in others, the empty  (``$\!(-1)\!$--dimensional'')
simplex is excluded, so that complexes are only assumed
closed under taking {\em nonempty} subsets of their members.
I quoted the former definition in \S\ref{S2.Delta} above.
But the choice makes only a formal difference to the study
of abstract simplicial complexes $D,$ and none to their geometric
realizations $\Delta(D),$ the subject of this note.

\subsection{Dependence among sections of this note}\label{S2.dependence}
All sections of this paper assume the material of \S\ref{S.constr}.
The remaining major sections,
\S\S\ref{S.variants}-\ref{S.stitching}, are largely independent
of one another, with the following exceptions:
\S\ref{S.more} assumes \S\ref{S.WT};
the construction of~\S\ref{S2.cl}
is used in~\S\S\ref{S2.thickenM3}-\ref{S2.thicken}
and in~\S\ref{S2.L-01}, and~\S\ref{S2.L-01}
also refers to the construction of~\S\ref{S2.Delta(P,S)}.

Within each of the main sections, later subsections
often depend on earlier ones.

\section{Generalizations and variants of the construction \texorpdfstring{$\Delta(-)$}{Delta}}\label{S.variants}

Before studying the construction $\Delta(-),$ let us digress and
record some ways it can be modified.

\subsection{Infinite \texorpdfstring{$L$}{L} or \texorpdfstring{$P\!$}{P}}\label{S2.inf}
We have been assuming that our given lattice $L$ is finite.
If we allow infinite $L,$ we have to decide, first,
from what topology on $[0,1]^L$ to induce a topology
on $\Delta(L)\subseteq [0,1]^L.$
The product topology will not, in general, make joins continuous.
For instance, suppose $L$ consists of
$0,$ $1,$ and infinitely many pairwise incomparable elements
$x_i$ $(i\in\mathbb{N})$ lying between them, so
that for $i\neq j$ we have $x_i\vee x_j=1$ and $x_i\wedge x_j=0.$
Then the sequence ${\downarrow}\,x_1,$
${\downarrow}\,x_2,$ $\dots$ in $\Delta(L),$
regarded as a sequence of $\![0,1]\!$--valued functions,
approaches ${\downarrow}\,0$
in the product topology on $[0,1]^L;$
but $\lim_{i\to\infty} ({\downarrow}\,x_0\vee{\downarrow}\,x_i) =
\lim_{i\to\infty} {\downarrow}\,1= {\downarrow}\,1,$ which is not
equal to ${\downarrow}\,x_0\vee (\lim_{i\to\infty} {\downarrow}\,x_i) =
{\downarrow}\,x_0\vee {\downarrow}\,0 = {\downarrow}\,x_0.$

If we use instead the topology on $[0,1]^L$ given
by the metric $d(\,f,g)=\sup_{x\in L} |f(x)-g(x)|,$ which for finite $L$
agrees with the product topology, the above problem goes away.
Indeed, note that if $d(\,f,f')<\varepsilon,$ then for all $t$
we have $f_{t+\varepsilon}\subseteq f'_t\subseteq f_{t-\varepsilon}$
(where we replace $t+\varepsilon$ by $1$ if $t+\varepsilon>1,$ and
$t-\varepsilon$ by $0$ if $t-\varepsilon<0).$
It is easy to deduce from this and the description of joins in
Theorem~\ref{T.main} that whenever
$d(\,f,f')<\varepsilon$ and $d(g,g')<\varepsilon,$ we get
$d(\,f\vee g,f'\vee g')<\varepsilon,$ and similarly for meets.

Using this topology,
if we assume that all chains in $L$ have finite length, then
the situation is quite similar to that of \S\ref{S.constr}.
The set of functions $f$ such that all the sets $f_t$ are principal
ideals will be a possibly infinite simplicial complex, but
still composed of finite-dimensional facets (maximal simplices).
Each simplex is still compact, though the whole space in general is not.
(For an example of noncompactness, let $L$ again
be the lattice whose elements are $0,$ $1$ and $x_i$ $(i\in\mathbb{N}).$
Then the elements
${\downarrow}\,x_i\in\Delta(L)$ have pairwise distance $1,$ hence
can have no convergent subsequence in our metric.
For another noncompact topology, see \cite[p.\,98, line~10]{JW}.)

If $L$ has infinite chains, there are more choices to be made.
We can again look at the set of $f\in[0,1]^L$ which
are finite convex linear combinations of characteristic functions
of chains of principal ideals; or at the larger
set of those $f$ such that every $f_t$
$(t\in\ho)$
is a principal ideal; these can be described as the {\em possibly
infinite}
convex linear combinations of chains of such characteristic functions.
The former space is still made up of finite-dimensional
simplices, though such
simplices will not in general be contained in maximal
simplices; the latter space is less like a simplicial complex.
In another direction, we can generalize principal ideals
to nonempty ideals $(\,=$ directed unions of principal ideals).
If we do so, we must again choose between using convex linear
combinations of finite chains of such ideals, and of arbitrary chains.

Most of the above considerations apply not only to lattices,
but also to $\!\vee\!$--semilattices, $\!\wedge\!$--semilattices,
and general posets $P.$
However, in the cases of $\!\wedge\!$--semilattices and posets, there
is no evident reason to reject the topology
induced by the product topology on $[0,1]^P$
in favor of the one based on
the metric $\sup_{x\in P} |f(x)-g(x)|.$
The meet operation is continuous under both topologies,
and the order-relation ${\leq},$ as a subset
of $\Delta(P)\times\Delta(P),$ is closed in both topologies.

\subsection{Constructions based on closure operators (etc.)}\label{S2.cl}
Returning for simplicity to finite lattices $L,$
what if we abstract the underlying set of $L,$ and its
family of principal ideals, to a
general finite set $X,$ and its closed sets under an arbitrary
closure operator $\r{cl}$ on $X$?

If we write $\Delta(X,\r{cl})$ for the set of functions
$f{:}\ X\to[0,1]$ such that for each
$t\in\ho,$
the set $f_t=\{x\in X\mid f(x)\geq t\}$ is closed under $\r{cl},$
we see that $\Delta(X,\r{cl})$ will again be a union of simplices,
indexed now by the chains of $\!\r{cl}\!$--closed subsets of $X;$
and that we can give it a lattice structure in the same
way we did for $\Delta(L).$
In fact, if we write $L_\r{cl}$ for the finite lattice
of $\!\r{cl}\!$--closed subsets of $X,$ the above
simplicial complex with lattice structure
will be isomorphic, as such, to $\Delta(L_\r{cl}).$
The isomorphism $\Delta(X,\r{cl})\to\Delta(L_\r{cl})$
can be described as sending each $f\in\Delta(X,\r{cl})$
to the function on $L_\r{cl}$ which takes each $A\in L_\r{cl}$
to the minimum value of $f$ on $A\subseteq X,$
or to $1$ if $A=\emptyset;$
its inverse sends each $f\in\Delta(L_\r{cl})$ to the function
taking each $x\in X$ to the maximum
value of $f$ on those members of $L_\r{cl}$ which (as subsets of $X)$
contain $x;$ equivalently, to $f(\r{cl}(\{x\})).$
So the construction $\Delta(X,\r{cl})$ does not give
{\em new} topological lattices; but at times we will find it
gives more convenient descriptions of the topological lattices we
are interested in than the construction~$\Delta(L).$

The above construction extends to the situation where we
are merely given a finite set $X$ and a family $P$ of subsets of $X.$
The functions $f{:}\ X\to[0,1]$ for which all of the sets $f_t$ belong
to $P$ forms a simplicial complex
isomorphic to $\Delta(P),$ where $P$ is regarded as
ordered by inclusion.
The isomorphism can be described as in the preceding
paragraph (with $\Delta(P)$ for $\Delta(L_\r{cl})),$ except that,
of the two equivalent
descriptions of one direction of the isomorphism, the
one using the sets $\r{cl}(\{x\})$ must be dropped.

If $L$ is a finite lattice, a convenient representation of
$L$ by closed subsets of a
set $X$ is gotten by taking for $X$ the set
of join-irreducible elements of $L$ other than the least
element (which we may regard as the join of the empty set), and for
$\r{cl}$ the operator taking every $S\subseteq X$ to
$X\cap{\downarrow}\,(\bigvee\,S).$

\subsection{Generalizing \texorpdfstring{$[0,1]\!$}{[0,1]}}\label{S2.M}
In a different direction, suppose we replace the interval
$[0,1]$ in our construction of $\Delta(L)$ with any
{\em complete lattice} $M.$
Thus, given a finite lattice $L,$ let $\Delta_M(L)$
denote the set of all set-maps $f{:}\ L\to M$ such that
for each $t\in M,$ the set $f_t=\{x\in L\mid f(x)\geq t\}$ is
a principal ideal of $L.$
(Our use of
$t\in\ho$
elsewhere in this note suggests that we say,
``for each $t\in M-\{0_M\}$''.
But as noted, this makes no difference since $L$ has a
greatest element.)
Partially ordering these functions $f$ by pointwise comparison,
we again see that meets of arbitrary families of elements can
be obtained as pointwise meets.
(This includes the meet of the empty family, which is
the constant function with value everywhere $1_M\in M.)$
It follows that every subset $S\subseteq\Delta_M(L)$ also
has a join, the meet of all members of $\Delta_M(L)$ that
majorize all members of $S.$
However, when $M$ is not a chain, we can no longer describe
$f\vee g$ by the condition that each $(\,f\vee g)_t$
be the join of $f_t$ and $g_t$ as ideals of $L.$
Indeed, let $L$ be a $\!2\!$--element lattice
$\{0_L,1_L\},$ and $f,\,g$ set-maps $L\to M$ whose values at
$0_L$ are both $1_M,$ and whose values at
$1_L$ are incomparable elements $p,\,q\in M.$
We see that $f, g\in\Delta_M(L);$ and
writing $t=p\vee q\in M,$ we see that the value of $f\vee g$
at $1_L$ will be $t,$ and deduce that $(\,f\vee g)_t={\downarrow}\,1_L,$
though $f_t\vee g_t={\downarrow}\,0_L\vee {\downarrow}\,0_L
={\downarrow}\,0_L.$

If we associate to every $f\in\Delta_M(L)$ the
set $\overline{f}=\{(x,y)\in L\times M\mid y\leq f(x)\},$
which clearly determines it, then we find that the
sets $\overline{f}\subseteq L\times M$ that we
get are characterizable by
a pair of conditions symmetric in $L$ and $M;$ namely,
they are those subsets such that
for each $x\in L,$ the set $\{y\in M\mid (x,y)\in\overline{f}\,\}$ is
a principal ideal of $M,$ and such that
for each $y\in M,$ the set $\{x\in L\mid (x,y)\in\overline{f}\,\}$ is
a principal ideal of $L.$
This is, in fact, a case of a known construction,
the {\em tensor product} $L\otimes M$ of lattices~\cite{GG+FW}.
So for $L$ a finite lattice, the lattice $\Delta(L)$
that is the main subject of this note
can be described as $L\otimes[0,1].$
If $L$ and $M$ are not both complete, their tensor
product may fail to be a lattice
(e.g., this happens if one is $M_3,$ and the other is a
free lattice on $3$ generators),
but in that situation, variants of the construction
have been described which do give lattices~\cite{GG+FW}.

The above construction $\Delta_M(L)$
may, of course, be generalized to allow $M$ and/or $L$
to be semilattices or posets rather than lattices.
In doing this, one must make choices on what to do when
one or both of these does not have a greatest element
(cf.\ \S\ref{S2.01}).

\subsection{A more up-down symmetric construction}\label{S2.sym}
Returning to the construction $\Delta(L),$
or more generally, $\Delta(P),$ note that our description of this
complex in terms of $\![0,1]\!$--valued functions
is very asymmetric with respect to the order relation:
$\Delta(P)$ consists of convex linear combinations of chains of
characteristic functions of principal {\em ideals}, rather than
principal {\em filters} (dual ideals).
A consequence is that for $L$ a lattice, though meets in $\Delta(L)$ are
pointwise infima, joins have a less trivial description.
In contrast, the standard description
of the abstract order complex $\Delta(P)$ in terms of chains of
elements of $P$ is completely symmetric.

We could, of course, {\em reverse}
the asymmetric feature of our construction of the spaces $\Delta(P)$ and
$\Delta(L),$ and use convex linear combinations of characteristic
functions of principal filters ${\uparrow}\,x$ instead of
principal ideals ${\downarrow}\,x.$
A minor difficulty is that larger elements generate
smaller filters; so the resulting lattice would
reverse the ordering of $L.$
However, we could fix this by interchanging $0$ and $1$ in
our characteristic functions; equivalently, by using characteristic
functions of complements of principal filters.
(The set of such complements has properties dual to those of a
closure system:
every subset of $L$ {\em contains} a largest member of this set.)

In fact, our original construction, and the above dual approach,
and the generalization in \S\ref{S2.cl} where the lattice of
principal ideals of $L$ was replaced by the lattice of
closed sets of any closure operator on a finite set,
can be subsumed under one construction.
The version of this construction for a general poset $P$ was
sketched in the last paragraph of~\S\ref{S2.cl}, so below, we
simply describe the lattice case.

\begin{lemma}\label{L.up+down}
Let $X$ be a finite set, and $L$ a family of subsets of $X$ which,
under the partial ordering by inclusion, forms a lattice; i.e.,
such that for every $a,b\in L$ there is a least member of
$L$ containing $a\cup b,$ denoted $a\vee b,$
and a greatest member of $L$
contained in $a\cap b,$ denoted $a\wedge b.$

Let $\Delta(X,L)$ be the set of functions $X\to [0,1]$
such that for every $t\in\ho,$
the set $f_t=\{y\in X\mid f(y)\geq t\}$ belongs to $L.$

Then $\Delta(X,L)$ is a simplicial complex with a
lattice structure, isomorphic as such to the object $\Delta(L)$
constructed in Theorem~\ref{T.main} from the abstract lattice $L.$
\end{lemma}

\begin{proof}[Sketch of proof]
It is easy to see that every element of $\Delta(L)$
has a unique representation as a convex linear combination
with nonzero coefficients
of the characteristic functions of a chain of members of $L,$
and that this gives an isomorphism as ordered simplicial complexes
with $\Delta(L).$

Let us describe explicitly the lattice operations of $\Delta(X,L).$
Given $f,\,g\in\Delta(L)$ the function
$f\wedge g\in\Delta(L)$ can be determined
by specifying the sets $(\,f\wedge g)_t$ for all $t.$
If we let $(\,f\wedge g)_t\subseteq X$ be the meet
$f_t\wedge g_t$ in $L$ of the sets $f_t,\,g_t,$
it is straightforward to verify that this gives
an element $f\wedge g\in\Delta(X,L),$ which
will be a greatest lower bound of $f$ and $g.$

Since our hypotheses are up-down symmetric,
the construction of $f\vee g$ is analogous.
(It is, in fact, the construction we used in Theorem~\ref{T.main}.)
\end{proof}

\Needspace{3\baselineskip}
\section{Two questions of Walter Taylor, an answer to one of them,\\and some approaches to the other}\label{S.WT}

\subsection{Modular but not distributive}\label{S2.M3}
As mentioned earlier, this paper was
motivated by the study by Walter Taylor \cite{WT} of topological spaces
admitting algebraic structures of various sorts.
One question he asked in an earlier version of \cite{WT} was whether
there exist topological spaces -- in particular, ``nice'' ones
such as simplicial complexes -- which admit structures
of modular lattice, but not of distributive lattice.

Consider the topological lattice
$\Delta(M_3),$ where $M_3$ is the lattice
\raisebox{2pt}[12pt][8pt]{ % M_3
\begin{picture}(15,20)
\dotline{7,-8}{7,0}
\dotline{7,0}{7,8}
\dotline{7,8}{0,0}
\dotline{0,0}{7,-8}
\dotline{7,8}{14,0}
\dotline{14,0}{7,-8}
\end{picture}}\,.

Since $M_3$ is modular but not distributive, the same is
true of $\Delta(M_3).$
In early drafts of this note, I could not go farther than
to suggest that the underlying
topological space of that object might not admit a distributive
lattice structure, in which case it would answer the above question.
Then Walter Taylor discovered that the
nonexistence of a distributive lattice structure
on that space can be proved by
combining two results in the literature.
With his permission, I give the argument here.
We first recall

\begin{definition}\label{D.breadth}
A finite nonempty subset of a lattice $L$ is said to be
{\em meet-irredundant}
if the meet of that set is not equal to the meet of any proper subset.
\textup{(}We shall regard the greatest element of $L,$ it this exists,
as the meet of the empty set.
Thus a singleton $\{x\}\subseteq L$ is meet-irredundant
if and only if $x$ is not a greatest element of $L.)$

The {\em breadth} $\r{br}(L)$ of $L$ is the supremum of the
cardinalities of all its meet-irredundant subsets, a natural number
or $+\infty.$
\textup{(}By the above observation on singleton families, the
$\!1\!$--element lattice has breadth $0.$
It is not hard to see that the lattices
of breadth $1$ are the chains with more than one element.\textup{)}
\end{definition}

The above definition is very meet-join asymmetric;
but it is easily shown to be equivalent to its dual.
In \cite[\S1.7]{GG.Lat}, a symmetric definition of the
breadth of an arbitrary partially ordered set is given,
which for lattices is equivalent to the above definition.

The next result refers to the ``inductive dimension'' of a
topological space.
This is a topological invariant about which we only need to know
that on simplicial complexes, it is equal to the usual dimension.

\begin{proposition}[Walter Taylor, personal communication]\label{P.WT}
Let $X$ be a compact connected metrizable topological space of
inductive dimension $n,$ such that the
set of points at which $X$ has dimension n has nonempty interior.
Then if $X$ admits a structure of distributive topological
lattice, $X$ is embeddable, as a topological space, in~$\mathbb{R}^n.$

Hence the underlying topological space of $\Delta(M_3)$
does not admit a structure of distributive lattice.
\end{proposition}

\begin{proof}
Let $L$ be a distributive topological lattice whose underlying space
has the properties of the first sentence.
Choe~\cite{Choe} shows that $\r{br}(L)=n,$ while Baker and
Stralka~\cite[Corollary~3.5]{B+S} show that
if $\r{br}(L)=n,$ then $L$ embeds homeomorphically
in an $\!n\!$--cell, hence in $\mathbb{R}^n.$
This gives the first assertion.
(The results cited from \cite{Choe} and~\cite{B+S} each
have a somewhat more general hypothesis; the hypothesis stated
above is roughly their intersection.)

On the other hand, the underlying space
of $\Delta(M_3),$ a $\!2\!$--dimensional simplicial complex
consisting of three triangles
meeting at a common edge, is not embeddable in $\mathbb{R}^2.$
This is intuitively clear, but here is a formal proof.
Call the three triangles in question $A,$ $B$ and $C.$
Then $A\cup B$ is homeomorphic to the closed unit
disk; let us call the subset thereof corresponding to
the open unit disc $U.$
Since $U$ is homeomorphic to an open subset of $\mathbb{R}^2,$
Invariance of Domain \cite[Theorem~6.36]{JR} tells us that
if $\Delta(M_3)$ were embedded in $\mathbb{R}^2,$ then
$U$ would be mapped to an open subset of $\mathbb{R}^2;$
so $U$ would have to be open in $\Delta(M_3).$
But $U$ is not open, since writing $E$ for
the edge joining $A,$ $B$ and $C,$ the set
$C-E,$ which is disjoint from $U,$
has in its closure the interior of $E,$ which is in $U.$

Hence by the first assertion of the proposition, $\Delta(M_3)$
cannot admit a structure of distributive topological lattice.
\end{proof}

The object $\Delta(M_3)$ has made earlier appearances in
the literature; see the picture on p.\,4 of~\cite{G+S}, and
the paragraph preceding that picture.

We remark that another way to prove that $\Delta(M_3)$ can't be embedded
in $\mathbb{R}^2$ is to show that it contains a non-planar graph.
If we give the vertices of $\Delta(M_3)$ the names $a,\,b,\,c,\,0,\,1,$
with the last two being the least and greatest elements,
and let $1/2$ be the midpoint
of the edge $E$ connecting $0$ and $1,$ then we see that
$\Delta(M_3)$ contains all the line-segments connecting
points of $\{a,b,c\}$ to points of $\{0,1/2,1\},$
giving a copy of the complete
bipartite graph $K_{3,3},$ which is non-planar
\cite[Corollary~11.1, cf.~Theorem~11.13]{FH}.
A more elaborate construction, shown to me by Taylor, embeds
the non-planar graph $K_5$ in $\Delta(M_3).$
For some very general results which include the embeddability of $K_5$
in $\Delta(M_3),$ see \cite{wiki_book}.

Before Taylor supplied Proposition~\ref{P.WT}, my thought on
how one might prove nonexistence of a distributive lattice structure on
the underlying space of $\Delta(M_3)$ centered
on the following question, to which I still do not know the answer.

\begin{question}\label{Q.3parts}
If $X$ is a connected finite simplicial complex having a subcomplex $Y$
such that the subspace $X-Y$ has at least $3$
connected components, each having all of $Y$ in its closure, can $X$
admit a structure of distributive topological lattice?
\end{question}

We end this section with a result of the same nature as the
one quoted from~\cite{Choe} in the proof of Proposition~\ref{P.WT},
but for not-necessarily-distributive lattices.
We could not have used it in that proof, since
it is restricted to lattices
of the form $\Delta(L),$ and we were interested in arbitrary
possible topological lattice structures on $\Delta(M_3).$
But the parallelism with Choe's result is interesting.

\begin{lemma}\label{L.breadth}
If $L$ is a finite lattice, then the breadth of the
lattice $\Delta(L)$ is equal to the supremum of the lengths
of all chains in $L$ \textup{(}where by the length of a finite
chain we understand one less than the number of elements in the
chain\textup{)}.
In other words \textup{(}in view of the structure of $\Delta(L)),$
$\r{br}(\Delta(L))=\r{dim}(\Delta(L)).$
\end{lemma}

\begin{proof}
For the easy direction, suppose that $0=x_0<\dots<x_n$ is a
chain in $L.$
Choose any real numbers $1>r_1>\dots>r_n\geq 0,$ and
for $i=1,\dots,n,$ let $f_i\in\Delta(L)$ be defined to have the
value $1$ on ${\downarrow}\,x_{i-1},$ and $r_i$ elsewhere.
Then we see that the value of $\bigwedge f_i$ at each $x_i$ is $r_i;$
but that if we omit any $f_i$ from this meet, the result
has the larger value $r_{i-1}$ at $x_i;$ so the meet is irredundant,
so $\Delta(L)$ has breadth at least~$n.$

For the converse, suppose $f_1,\dots,f_n\in\Delta(L)$ are
distinct elements forming a meet-irredundant set, and let us
show that $L$ has a chain of length $n.$

For each $i,$ we have $f_i\not\geq \bigwedge_{j\neq i}\,f_j,$
so for each $i$ there must exist $x_i\in L$ such
that $f_i(x_i)<f_j(x_i)$ for all $j\neq i.$
Let us choose such an $x_i$ for each $i,$ and reindex our
elements, if necessary, so that
\begin{equation}\begin{minipage}[c]{25pc}\label{d.f1x1}
$f_1(x_1)\ \geq\ \dots\ \geq\ f_n(x_n).$
\end{minipage}\end{equation}
Then I claim that in $L,$
\begin{equation}\begin{minipage}[c]{25pc}\label{d.<<}
$0\ <\ x_1\ <\ x_1\vee x_2\ <\ \dots\ <\ x_1\vee\dots\vee x_n.$
\end{minipage}\end{equation}
Since the relation $\leq$ holds between successive terms
of~\eqref{d.<<}, it suffices to show that for each $i$
we have $\bigvee_{j<i} x_j\neq\bigvee_{j\leq i} x_j.$
We shall do this by showing that $f_i$ has distinct values
at these two elements of $L.$
By~\eqref{d.vee=wedge},
$f_i$ turns joins in $L$ into infima in $[0,1],$ so what we must show
is that $\bigwedge_{j<i}\,f_i(x_j)\neq\bigwedge_{j\leq i}\,f_i(x_j);$
in other words,
that each term $f_i(x_j)$ $(j<i)$ is larger than $f_i(x_i).$
And indeed, $f_i(x_j)>f_j(x_j)\geq f_i(x_i),$ the
first inequality by choice of $x_j,$ the second by~\eqref{d.f1x1}.
Thus, $L$ has a chain~\eqref{d.<<} of length $n,$ as desired.
\end{proof}

This suggests

\begin{question}\label{Q.breadth}
Can the equality $\r{br}(M)=\r{dim}(M)$ be proved for wider classes
of topological lattices $M$ than
\textup{(i)}~those of the form $\Delta(L),$
as shown in Lemma~\ref{L.breadth}, and
\textup{(ii)}~the distributive topological lattices for which
the same is shown in \cite{Choe}?
\textup{(}E.g., for all topological lattices whose underlying spaces
are connected simplicial complexes, or still more generally,
all connected compact Hausdorff spaces?\textup{)}

Can one at least obtain one or the other inequality between these
invariants?
\end{question}

\subsection{What about an example with no modular lattice structure?}\label{S2.N5}
Recall that the modular
lattice $M_3$ is one of two $\!5\!$--element non-distributive
lattices, the other being the nonmodular lattice $N_5=$
\raisebox{2pt}[12pt][8pt]{ % N_5
\begin{picture}(15,20)
\dotline{6,-8}{0,-3.5}
\dotline{0,-3.5}{0,3.5}
\dotline{0,3.5}{6,8}
\dotline{6,8}{12,0}
\dotline{12,0}{6,-8}
\end{picture}}.
Thus, we might hope that $\Delta(N_5)$ would give an
example answering the other half of
Taylor's question on lattice identities:
whether a topological space can admit a structure of
lattice, but not of modular lattice \cite[Question~9.4.7]{WT}.
As a simplicial complex, $\Delta(N_5)$ consists of a tetrahedron
with a triangle attached along one edge.
Now it is certainly true that a simplicial complex $\Delta(L),$ for
$L$ a modular lattice, cannot have this form; for in a finite
modular lattice $L,$ all maximal chains have the same length
\cite[Theorem~374]{GG.Lat}, hence in its order complex,
all facets (maximal simplices) have the same dimension.
Nevertheless, it is easy to construct a topological distributive
(and hence modular) lattice {\em not} of the form $\Delta(L)$
whose underlying space is homeomorphic to $\Delta(N_5).$
Simply take the $\!2\!$--dimensional lattice $[0,1]^2,$ ``glue''
a copy of $[0,1]$ to its top (identifying the top element of
$[0,1]^2$ with the bottom element of $[0,1]),$ then take the
direct product of the resulting lattice with the lattice $[0,1].$
Geometrically, the result is a cube with a square attached along one
edge; and since a cube is homeomorphic to a tetrahedron, by a
homeomorphism that can be made to carry one edge
of the cube to an edge of the tetrahedron, and a square
is similarly homeomorphic to a triangle, the above example
is indeed homeomorphic to $\Delta(N_5).$

On the other hand, if one started with a lattice $L$ such as
\raisebox{2pt}[12pt][8pt]{ % 2+1+1
\begin{picture}(15,20)
\dotline{6,-8}{0,-3.6}
\dotline{0,-3.6}{0,3.6}
\dotline{0,3.6}{6,8}
\dotline{6,8}{13,0}
\dotline{13,0}{6,-8}
\dotline{6,-8}{6,0}
\dotline{6,0}{6,8}
\end{picture}}\,,
or, if need be,
\raisebox{2pt}[12pt][8pt]{ % 2+1+2
\begin{picture}(15,20)
\dotline{6,-8}{0,-3.6}
\dotline{0,-3.6}{0,3.6}
\dotline{0,3.6}{6,8}
\dotline{6,8}{12,3.6}
\dotline{12,3.6}{12,-3.6}
\dotline{12,-3.6}{6,-8}
\dotline{6,-8}{6,0}
\dotline{6,0}{6,8}
\end{picture}}
or
\raisebox{2pt}[12pt][8pt]{ % 2+2+2
\begin{picture}(15,20)
\dotline{6,-8}{0,-3.0}
\dotline{0,-3.0}{0,3.0}
\dotline{0,3.0}{6,8}
\dotline{6,8}{12,3.0}
\dotline{12,3.0}{12,-3.0}
\dotline{12,-3.0}{6,-8}
\dotline{6,-8}{6,-3.0}
\dotline{6,-3.0}{6,3.0}
\dotline{6,3.0}{6,8}
\end{picture}}\,,
we can hope that the underlying topological space
of $\Delta(L)$ will not admit any modular lattice structure.
The intuition is that although when one glues one topological lattice
on top of another, as we did with a square and a line-segment above,
the common connecting sublattice can have much lower dimension than the
two parts we are connecting, there is no evident construction
that preserves distributivity {\em or} modularity and
unites {\em more than two} parts along a common connection $Y$ having
codimension greater than~$1.$
So we ask,

\begin{question}\label{Q.3parts&codim}
Suppose $X$ and $Y$ are as in Question~\ref{Q.3parts},
and moreover suppose that at least one \textup{(}at least two?
at least three?\textup{)}\ of the connected components
of $X-Y$ having $Y$ in their
closures have dimension exceeding that of $Y$ by at least two.
Can $X$ admit a structure of {\em modular} topological lattice?
\end{question}

The order complexes of the three nonmodular lattices
shown above have the properties described in the above question,
so a positive answer to any of the three versions
of that question would answer this part of Taylor's question.
We will look at another example that might do so in~\S\ref{S2.whiskers}.

\subsection{A construction suggested by the glued example}\label{S2.Delta(P,S)}
The distributive lattice constructed in the
preceding section by gluing $[0,1]$ onto
$[0,1]^2$ exhibits a property that we noted
could not occur in the order complex $\Delta(L)$ of a modular
lattice $L.$
Thus, not every simplicial complex
with lattice structure arises as a $\Delta(L).$
Can we modify the construction $\Delta(L)$ to embrace this example?
In a somewhat ad hoc way, we can.

\begin{proposition}\label{P.pairs}
If $P$ is a finite partially ordered set, and $S\subseteq P\times P$ a
family of pairs $(x,y)$ each satisfying $x\leq y,$ let $\Delta(P,S)$ be
the subspace of $\Delta(P)$ consisting of functions
$f\in\Delta(P)$ with the property that for all
$(x,y)\in S,$ if $f(x)<1$ then
$f(y)=0$ \textup{(}equivalently, such that $f(x)$ and
$f(y)$ do not both lie strictly between $0$ and~$1).$

Then\vspace{.2em}

\textup{(i)}\ \ $\Delta(P,S)$ is a subcomplex of $\Delta(P)$
with the same set of vertices.

\textup{(ii)}\ \ If $P$ has meets of all pairs
of elements, then the meet operation of $\Delta(P)$ carries
$\Delta(P,S)$ into itself.

\textup{(iii)} If $P$ has joins of all pairs of
elements, and if the second coordinates $y$ of the pairs $(x,y)\in S$
are all join-prime \textup{(}definition recalled
below\textup{)}, then the join operation of
$\Delta(P)$ carries $\Delta(P,S)$ into itself.\vspace{.2em}

Thus, if $P$ is a lattice $L,$
and the second coordinates of all members of $S$ are join-prime
in $L,$ then $\Delta(L,S)$ is a sublattice of $\Delta(L),$
containing the sublattice of vertices.
\end{proposition}

\begin{proof}
Given $f\in\Delta(P),$ if we write $f$ as a convex linear
combination with nonzero coefficients of the characteristic functions
of the principal ideals generated by elements
$z_0<\dots<z_n,$ then we find that
$f$ will belong to $\Delta(P,S)$ if and only if for each $(x,y)\in S,$
either $z_0\geq x$ (which is equivalent to
$f(x)=1)$ or $z_n\not\geq y$ (equivalent to $f(y)=0).$
The class of nonempty chains $z_0<\dots<z_n$ satisfying this
condition for a given $(x,y)\in S$
is clearly closed under taking nonempty subchains, hence the
points $f$ associated with chains having this property
for all $(x,y)\in S$ form a subcomplex
of $\Delta(P).$

At every {\em vertex} of $\Delta(P),$ all coordinates are $0$
or $1,$ so the condition that for all $(x,y)\in S,$ $f(x)$ and
$f(y)$ do not both lie strictly between $0$ and $1$
is certainly satisfied, completing the proof of~(i).

If $P$ has meets, as assumed in (ii), so that $\Delta(P)$ is
closed under pointwise infima, I claim that
$\Delta(P,S)$ is also closed under these infima.
Indeed, for $f,g\in \Delta(P,S)$ and $(x,y)\in S,$
either $f(x)$ and $g(x)$ are both $1,$ in which
case the infimum $f\wedge g$
certainly satisfies the required condition, or one of them, say $f(x),$
is $<1,$ in which case $f(y)=0,$ and the infimum again satisfies
our condition.

Finally, suppose that $P$ has pairwise joins, and that for
each $(x,y)\in S,$ the element $y$ is join-prime,
meaning that $w\vee z\geq y$ only if $w\geq y$ or $z\geq y.$
Again let $f,g\in\Delta(P,S)$ and $(x,y)\in S.$
If $(\,f\vee g)(y)=t>0,$ then by the construction
of joins in the semilattice
$\Delta(P),$ we must have $f(y_1)\geq t$ and $g(y_2)\geq t$
for some $y_1,\,y_2$ with join $\geq y.$
Since $y$ is join-prime, one of $y_1,\,y_2$ must be $\geq y;$
say $y_1\geq y.$
Then $f(y)\geq f(y_1)\geq t>0,$ so $f(x)=1,$ so $(\,f\vee g)(x)=1,$
as required to show $f\vee g\in\Delta(P,S).$

The final assertion of the proposition is now clear.
\end{proof}

To see that the lattice obtained by gluing
a copy of $[0,1]$ to the top of a copy of $[0,1]^2$
is an example of this construction,
let $L$ be the $\!5\!$--element lattice
\raisebox{2pt}[12pt][8pt]{ % | glued on <>
\begin{picture}(20,20)
\dotline{10,-10}{5,-4}
\dotline{5,-4}{10,2}
\dotline{10,2}{15,-4}
\dotline{15,-4}{10,-10}
\dotline{10,10}{10,2}
\put(-3,-5){${\scriptstyle a}$}
\put(15,-5){${\scriptstyle b}$}
\put(10,2){${\scriptstyle x}$}
\put(10,8){${\scriptstyle y}$}
\end{picture}}\,,
and $S$ the singleton $\{(x,y)\}.$
We see that the elements of $\Delta(L)$ with $f(y)=0$ form a sublattice
isomorphic to $[0,1]^2,$ each element being determined by its
values at $a$ and $b,$ while the elements with $f(x)=1$
form a sublattice isomorphic to $[0,1],$
each being determined by its value at $y.$
These sublattices together comprise $\Delta(L,S),$ and the
greatest element of the former is the least
element of the latter, giving the asserted structure.

One can, in turn, get from this the cube-with-a-square-attached,
as $\Delta(L{\times}\{0,1\},\,S')$ where $S'=\{((x,0),(y,0))\}.$

Incidentally, the need for the condition in Question~\ref{Q.3parts}
that each of the connected components of $X-Y$ considered have
$Y$ in its closure is illustrated by the
line-glued-on-top-of-a-square lattice $\Delta(L,S)$ of the
next-to-last paragraph above.
Though $\Delta(L,S)$ is distributive, if we denote by $0$
the least element of $L,$ then
removing the $\!1\!$--simplex $Y=\Delta(\{0,x\}),$
i.e., the vertical diagonal of the square, breaks
$\Delta(L,S)$ into three pieces, whose closures are respectively
the triangles $\Delta(\{0,a,x\})$ and
$\Delta(\{0,b,x\})$ and the line-segment $\Delta(\{x,y\});$
but the last of these closures does not contain all of~$Y.$

We remark that in the context of the first sentence
of Proposition~\ref{P.pairs}, for any
four elements $x\leq x'\leq y'\leq y$ of $P,$ the condition on points
of $\Delta(P)$ determined by the pair $(x',y')$ implies
the condition determined by $(x,y);$ so if $(x',y')$ lies
in $S,$ it makes no difference whether $(x,y)$ does.
Hence we can restrict $S$ in that construction to
be an antichain in the poset of pairs of comparable elements, ordered by
inclusion of intervals.

\subsection{Lattices with whiskers}\label{S2.whiskers}

In \S\ref{S2.N5} we noted some variants of the topological
lattice $\Delta(N_5)$ which might have the property that their
underlying spaces would not admit a modular lattice structure.

There is another very simple example that might have this
property: the topological lattice gotten by gluing a copy of the unit
interval to the top of $\Delta(N_5).$
(I.e., the union of
$\Delta(N_5)$ with $[0,1],$ where the greatest element of the former
lattice is identified with the least element of the latter.
In the notation for gluing of \cite[end of \S{}I.1.6]{GG.Lat}, this
is $\Delta(N_5)\,\dot{+}\,[0,1];$ in the notation of
\S\ref{S2.Delta(P,S)} above, it is $\Delta(L,S),$
where $L$ is the result of adjoining to $N_5$ a new element $1^+$
above its existing top element $1,$ and $S$ is the
singleton $\{(1,1^+)\}.)$

Note that if we take the {\em distributive} topological lattice
which we saw in~\S\ref{S2.N5} was {\em homeomorphic} to $\Delta(N_5),$
and similarly glue a copy of $[0,1]$ to its top,
the result is a distributive lattice, but is {\em not} homeomorphic
to the space just described: the top of our distributive lattice
lies on the outer boundary of its $\!2\!$--dimensional part,
while the top of $\Delta(N_5)$
is one end of the edge along which its $\!2\!$--dimensional
and $\!3\!$--dimensional parts meet; so gluing $[0,1]$ to the
respective tops of these lattices gives nonhomeomorphic
topological spaces.

We shall see below that
any topological lattice structure on the underlying
space of our above extension of $\Delta(N_5)$ has the property
that the homeomorphic copy of $\Delta(N_5)$
is again a sublattice, and that
the point to which we have attached $[0,1]$ is either the
greatest or the least element of that lattice.
This might prove useful in showing that no such lattice
structure can be modular.

The key tool will be

\begin{proposition}\label{P.glue}
Suppose a topological lattice $L$ is connected, and has a point
$p$ such that $L-\{p\}$ is disconnected.

Then $L-\{p\}$ has precisely two connected components,
$({\downarrow}\,p)-\{p\}$ and $({\uparrow}\,p)-\{p\}.$
\textup{(}So in the notation of \cite[\S{}I.1.6]{GG.Lat},
$L={\downarrow}\,p\ \dot{+}\ {\uparrow}\,p.)$
\end{proposition}

\begin{proof}
Let us first prove that any two elements of $L$ neither of
which lies in ${\downarrow}\,p$ must
lie in the same connected component of $L-\{p\}.$
Indeed, if $x_0$ and $x_1$ are such elements, then
the map $f{:}\ L\to L$ defined by $f(y)=x_0\vee y$
has range in $L-\{p\},$ and has connected image since
it has connected domain; so $f(x_0)=x_0$ and $f(x_1)=x_0\vee x_1$
lie in the same connected component of $L-\{p\}.$
Similarly, $x_1$ and $x_0\vee x_1$ lie in
the same connected component of that space; so
$x_0$ and $x_1$ lie in the same connected component, as claimed.

The dual argument shows that any two elements not in ${\uparrow}\,p$
lie in the same connected component.

Since every element of $L-\{p\}$ is either ${\notin}\,{\downarrow}\,p$
or ${\notin}\,{\uparrow}\,p,$ that space has at most two
connected components.
Moreover, if any element were both ${\notin}\,{\downarrow}\,p$ and
${\notin}\,{\uparrow}\,p,$
those connected components would be the same, contradicting
our assumption of disconnectedness; so every element lies in one or
the other, and we see that the two connected components
are just $({\downarrow}\,p)-\{p\}$ and $({\uparrow}\,p)-\{p\}.$
\end{proof}

(The fact that in the context of the above result, $L-\{p\}$
cannot have more than two connected components is proved by
Taylor \cite[third through sixth paragraphs of \S7.2.1]{WT},
generalizing a result asserted by A.\,D.\,Wallace \cite{ADW}.)

The above result strongly restricts
lattice structures on topological spaces $S$ having a point $p$ such
that $S-\{p\}$ has precisely two connected components, $S_1$ and $S_2.$
We shall now apply it in the case where one of $S_0\cup\{p\}$ and
$S_1\cup\{p\}$ is a space we are interested in, while the other
is the simplest nontrivial connected topological space, $[0,1].$

\begin{corollary}\label{C.whiskers}
Let $X$ be a nontrivial connected topological space, $p$ a point of $X,$
and $X'$ the space obtained by attaching a
``whisker'' $W$ to $X$ at $p;$ that is, by
taking the disjoint union of $X$
with a copy $W$ of the unit interval, identifying $p\in X$
with one endpoint of $W,$ and giving the resulting space
the natural topology.

Then structures of topological lattice on $X'$ correspond
bijectively to structures of topological lattice on $X$
under which $p$ is either a greatest or a least element.
Namely, if $p$ is the greatest element of $X$
under a topological lattice structure, then in the corresponding
structure on $X',$ the ``whisker'' $W$ is a lattice-theoretic copy
of the unit interval, glued at its bottom to the top of that lattice,
while if $p$ is the least element of our lattice structure
on $X,$ then $W$ is a copy
of $[0,1]$ glued at its top to the bottom of that lattice.

Likewise, suppose $p,$ $q$ are points of $X,$ and let
$X''$ be the space obtained by attaching to $X$ one whisker, $W,$
at $p,$ and another whisker, $W',$ at $q.$
Then $X''$ cannot admit a structure of topological lattice
unless $p\neq q,$ and in that case, such structures correspond
to structures of topological lattice on $X$ under
which one of $p,\ q$ is the greatest element and the other
is the least element.
In the corresponding structures
on $X'',$ the whiskers $W$ and $W'$ are each order-isomorphic
to $[0,1],$ glued to $X,$ as described above, at $p$ and $q.$
\end{corollary}

\begin{proof}[Sketch of proof]
If we are given a point $p\in X,$ it is clear that
for any structure of topological lattice on $X$ under which $p$
is the greatest element, the result of ``gluing'' a copy of $[0,1]$
to the top of $X$ will give a topological lattice homeomorphic to $X';$
and we have the obvious corresponding
construction if $p$ is the least element of our lattice structure.
Conversely, if we are given a lattice structure on $X',$ then,
noting that the space $X'-\{p\}$ is disconnected, namely, that
it is the union of its open subsets
$S_0=X-\{p\}$ and $S_1=W-\{p\},$ Proposition~\ref{P.glue}
tells us that one of $X,$ $W$ is precisely
${\downarrow}\,p,$ and the other is ${\uparrow}\,p.$
By the symmetry of the result we want to prove, it will suffice
to consider the case where $X={\downarrow}\,p,$ $W={\uparrow}\,p.$

Clearly, this makes $X$ and $W$ sublattices of $X',$ such that
$X'$ is the lattice formed by gluing $W$ on top of $X.$
It remains only to say why the lattice structure on $W$ must agree
with the usual lattice structure on $[0,1].$
It is probably known that, up to antiisomorphism,  that is the only
topological lattice structure on that space; in any case, that
description of the lattice structure of $W$
is not hard to deduce here by noting that for any interior point
$p'$ of $W,$ deletion of $p'$ also disconnects $X',$ so
Proposition~\ref{P.glue} can be applied again to show that all
the elements of $W$ on the ``$X$'' side of $p'$ are $<p',$
and those on the other side are $>p'.$

To get the assertions of the
final paragraph, first note that if $p=q,$ then $X''$
falls into at least three pieces on deleting $p,$ hence by
Proposition~\ref{P.glue} it cannot admit a structure of topological
lattice.
Assuming that $p\neq q$ and that we are given a topological
lattice structure on $X'',$
the desired result can be obtained by applying twice the one-whisker
result proved above: first to $X'',$ viewed as obtained by adjoining
the whisker $W'$ to the space $X',$ and then to $X',$
viewed as obtained by adjoining the whisker $W$ to~$X.$
(In the latter case, we apply our result to
the lattice structure on $X'$ that the first step
shows us it inherits from $X'',$ as the sublattice
${\uparrow}\,q$ or~${\downarrow}\,q.)$
\end{proof}

It follows from the first statement of the above corollary
that if the underlying
topological space of $\Delta(N_5)$ admits no modular
topological lattice structure having the {\em same greatest
element} as that of $\Delta(N_5)$ itself, then
the result of gluing a whisker to the top of $\Delta(N_5)$
will be a nonmodular topological lattice whose underlying space
cannot be given {\em any} structure of modular topological lattice,
and so will give our desired example.
(Here I am using the fact that if the underlying space
of $\Delta(N_5)$ admits no modular
lattice structure having the same greatest element as its standard
structure, then by symmetry, it also admits no modular
lattice structure having for least element the
greatest element of that structure.)

We can hedge our bets by noting that
the second assertion of the corollary tells us
that unless the underlying space of $\Delta(N_5)$ admits a
modular lattice structure having {\em both} the same greatest element
{\em and} the same least element as in $\Delta(N_5),$ then the
result of attaching one whisker to the top of $\Delta(N_5)$ and
another to the bottom will be an example of the desired sort.

One can also use Proposition~\ref{P.glue} to show that the
condition which Proposition~\ref{P.WT} proves is necessary
for a simplicial complex to admit a {\em distributive} lattice structure
is not sufficient, even when $X$ does admit some lattice structure.
Namely, if we glue together the topological lattices $\Delta(M_3)$
and $[0,1]^3,$ identifying the top element of one with
the bottom element of the other, the result will be a topological
lattice whose underlying space has dimension~$3$
and is clearly embeddable in $\mathbb{R}^3;$
but if it admitted a distributive topological lattice
structure, Proposition~\ref{P.glue} shows that its
$\Delta(M_3)$ part would also, which we have seen it does not.
\vspace{.2em}

The method of proof of Proposition~\ref{P.glue} also yields information,
though not as strong, about topological {\em semilattices}.
If $P$ is a connected $\!\vee\!$--semilattice which becomes
disconnected on deleting a point $p,$ one sees as in that proof
that any two elements ${\notin}\,{\downarrow}\,p$ lie in the same
connected component of $P-\{p\};$ in other words, all but
at most one of its connected components lie in ${\downarrow}\,p.$

Examples of such semilattices with many connected components
in ${\downarrow}\,p$ can be constructed by
starting with an arbitrary family of connected topological
$\!\vee\!$--semilattices $P_i$ $(i\in I),$ each having a greatest
element $p_i,$ and possibly one additional connected
topological $\!\vee\!$--semilattice $P_0$ $(0\notin I)$
with an arbitrary point $p_0$ chosen, identifying all the points $p_i$
(for $i\in I$ or $I\cup\{0\}$ depending on whether we have a $P_0),$
and giving the resulting set $P$ the topology and order-structure
determined in the obvious ways by those of the given structures.
(In particular, the join of elements $x_i\in P_i,$
$x_j\in P_j$ where $i\neq j$ and $i,j\in I$ is $p;$
while if there is a $P_0,$ the join of elements
$x_i\in P_i,$ $x_0\in P_0$ is $p\vee x_0\in P_0.)$

However, not all examples have this simple form.
For instance, consider the $\!\vee\!$--subsemilattice $P$ of $[0,1]^2$
consisting of the points $(x,x)$ and $(x,1)$ for all $x\in[0,1];$
and take $p=(1,1).$
Then $P-\{p\}$ has two connected components, but its
join operation does not have the form described above;
e.g., $(0,0)\vee(0,1)=(0,1),$ rather than being $p=(1,1).$

Returning to Proposition~\ref{P.glue} itself, this can be
generalized by replacing the singleton $\{p\}$ with
any convex sublattice $P\subseteq L$ (any sublattice
$P$ such that whenever
$x<y\in P,$ one has $[x,y]\subseteq P),$ such that
$L-P$ is disconnected.
Again, $L-P$ will have just two connected components.
The unions of these components with $P$ will be the sublattices
$\{x\in L\mid(\exists\,p\in P)\ x\leq p\}$
and $\{x\in L\mid(\exists\,p\in P)\ x\geq p\}.$
In the language of \cite[\S{}IV.2.1]{GG.Lat}, $L$ will be the lattice
obtained by gluing these two sublattices together over~$P.$

\section{More questions and examples}\label{S.more}

\subsection{A question on local distributivity}\label{S2.Q.loc_dist}
A notable property of lattices of the form $\Delta(L),$
and hence also of their sublattices $\Delta(L,S),$
is that each of their simplices is a distributive sublattice.
(Each simplex has the form $\Delta(C)$ for $C$ a chain in $L,$
and every chain $C$ in a lattice is a distributive sublattice;
hence so is the induced lattice $\Delta(C).)$
Thus, every interior point of every maximal simplex of such a
topological lattice has a neighborhood which is a distributive
sublattice.
This suggests

\begin{question}\label{Q.loc_dist}
Let $L$ be a topological lattice whose underlying set is
a finite simplicial complex \textup{(}or perhaps belongs to some
wider class, such as finite CW-complexes\textup{)}.
Suppose we call an element $p\in L$ a ``point of local distributivity''
if $p$ has a neighborhood which is a distributive sublattice of $L.$
Must the set of points of local distributivity \textup{(}an open
subset of $L)$ be dense in $L$?
\end{question}

\subsection{A short-lived variant of the above question, and another construction}\label{S2.thickenM3}
The observations that led to Question~\ref{Q.loc_dist}
also suggest a simpler question: whether every topological lattice
which is homeomorphic to an $\!n\!$--cell is distributive.

None of the examples we have seen contradicts this, but
a counterexample has long been known~\cite{DE}.
We will discuss that example in~\S\ref{S2.DE};
here I will describe a simpler one,
obtained by slightly ``thickening'' the lattice $\Delta(M_3).$
In developing this example it will be convenient to represent
$\Delta(M_3)$ as $\Delta(X,\r{cl})$ (notation as in~\S\ref{S2.cl}),
where $X=\{0,1,2\}$ and $\r{cl}$ is the
closure operator under which the singletons are closed, while the
closure of each $\!2\!$--element set is all of $X.$
Using this representation of $\Delta(M_3)$ will
have the effect that we will be ``thickening'' $\Delta(M_3)$
within the $\!3\!$--cube $[0,1]^3,$ rather than
having to work in the $\!5\!$--cube $[0,1]^L.$

To see what we will mean by ``thickening'',
note that for this $X$ and $\r{cl},$
$\Delta(X,\r{cl})$ consists of all functions
$f{:}\ \{0,1,2\}\to[0,1]$ such that if the value of $f$ at two of
these elements is $\geq t,$ then so is its value at the third.
This forces the  two lowest values of $f$ to be equal, so that
$\Delta(X,\r{cl})$ is $\!2\!$--dimensional.
Now let us fix a positive constant $c<1,$ and let $\Delta_c(X,\r{cl})$
denote the set of those functions
$f{:}\ \{0,1,2\}\to[0,1]$ such that if the value of $f$ at
two points is $\geq t,$ then its value at the third is $\geq t-c.$
It is not hard to see visually that this space is
homeomorphic to the $\!3\!$--ball;
one can give a formal proof starting with the observation that
every ray in $[0,1]^3$ through the point $(1/2,\,1/2,\,1/2)$ meets the
space in a nontrivial closed interval.
To see that it is a lattice under the ordering by componentwise
comparison, note that a least upper bound of any two elements
(e.g., $(1,0,0)$ and $(0,0,1))$ can be gotten by taking their
pointwise supremum and, if the value of one of the coordinates
thereof is
less than the minimum of the other two by more than $c,$ increasing
it to precisely that minimum minus $c.$
(So the least upper bound of $(1,0,0)$ and $(0,0,1)$ is $(1,1{-}c,1).)$
Greatest lower bounds are, as usual, calculated coordinatewise.

This lattice structure is still nondistributive, since
we find that
\begin{equation}\begin{minipage}[c]{25pc}\label{d.D_c_nondist}
$(0,1,0)\,\wedge\,((1,0,0)\vee(0,0,1))\ =\ (0,1{-}c,0),$ \\[.2em]
while \\[.2em]
$((0,1,0)\wedge(1,0,0))\,\vee\,((0,1,0)\wedge(0,0,1))\ =\ (0,0,0).$
\end{minipage}\end{equation}

Thus, it answers the suggested question in the negative.
\vspace{.2em}

Curiously, this lattice, though obtained by modifying a
lattice isomorphic to $\Delta(M_3),$ contains no copies of $M_3.$
Indeed, suppose $p_0,$ $p_1,$ $p_2\in\Delta_c(X,\r{cl})$ satisfy
\begin{equation}\begin{minipage}[c]{25pc}\label{d.p012}
$p_0\vee p_1= p_0\vee p_2=p_1\vee p_2,$ \quad
and \quad $p_0\wedge p_1= p_0\wedge p_2=p_1\wedge p_2.$
\end{minipage}\end{equation}
Let $a\in[0,1]$ be the greatest of
the $3\times 3=9$ coordinates of $p_0,$ $p_1$ and $p_2,$
and suppose without loss of generality that this
value $a$ is attained at the $\!0\!$--coordinate of $p_0.$
Then the $\!0\!$--coordinate of the element
$p_0\vee p_1= p_0\vee p_2=p_1\vee p_2$ must be $a,$
as can be seen from either of the first two expressions.
Hence $a$ is also the $\!0\!$--coordinate of the last expression,
$p_1\vee p_2,$ and, being the largest coordinate in any
of the $p_i,$ it cannot arise in that join
by the operation of increasing the smallest coordinate of
a $\!3\!$--tuple to the value of the second-largest coordinate minus $c$
(because $c>0);$ so it must be the $\!0\!$--coordinate of $p_1$ or $p_2.$
Assume without loss of generality that it is
the $\!0\!$--coordinate of $p_1.$
Then $a$ will also be the $\!0\!$--coordinate of
$p_0\wedge p_1= p_0\wedge p_2=p_1\wedge p_2$
(because it is the $\!0\!$--coordinate of the first of these);
moreover, since $a$ is the largest of the coordinates of the $p_i,$
the coordinate $a$ in these meets
cannot arise as the lesser of two distinct values; so
it must be the $\!0\!$--coordinate of {\em all} of $p_0,$ $p_1,$ $p_2.$
So these three elements have the same $\!0\!$--coordinate.
Next, let $b$ be largest of the $\!1\!$-- and
$\!2\!$--coordinates of $p_0,$ $p_1$ and $p_2,$ and assume
without loss of generality that it occurs as a $\!1\!$--coordinate.
Then the same argument shows that $b$ is the $\!1\!$--coordinate of all
three elements (since the second-largest coordinate also
cannot arise by the ``increase the smallest coordinate'' operation).
So $p_0,$ $p_1,$ $p_2$ agree in their first {\em two} coordinates.
Hence they form a chain in $\Delta_c(X,\r{cl});$ so the
relations~\eqref{d.p012} force $p_0=p_1=p_2,$ completing the proof that
$\Delta_c(X,\r{cl})$ contains no copy of~$M_3.$

Being nondistributive by~\eqref{d.D_c_nondist}, and containing
no $M_3,$ the lattice $\Delta_c(X,\r{cl})$ must contain an $N_5.$
And, indeed, if
$c\in(0,1/2]$
it is easy to verify that
the two elements $(1,0,0)<(1,c,0)$ belong to $\Delta(X,\r{cl}),$
and have equal joins with $(0,0,1),$
and also equal meets with that element; while if
$c\in[1/2,1)$
the same is true of the two elements $(1,0,0)<(1,1-c,0).$
In each case, this gives an explicit $\!5\!$--element
sublattice of the indicated form.
So though $\Delta(X,\r{cl})\cong\Delta(M_3)$ is, like $M_3,$
modular, the thickened lattice $\Delta_c(X,\r{cl})$ is not.

In view of the fact that $\Delta_c(X,\r{cl})$ contains
no copy of $M_3,$ the `short-lived question' with which
this section began might be revived in the weakened form:
Must a {\em modular} topological lattice that is homeomorphic
to an $\!n\!$--cell be distributive?
This is known to be true for $n\leq 3$ \cite{LA}~\cite{G+S}.

\subsection{The general thickening construction (and variants)}\label{S2.thicken}
Abstracting the construction of the preceding section, let us, for any
closure operator $\r{cl}$ on a finite set $X,$ and any
$c\in[0,1],$ define
\begin{equation}\begin{minipage}[c]{25pc}\label{d.D_c}
$\Delta_c(X,\r{cl})\ =\ \{\,f{:}\ X\to[0,1]\mid (\forall\,t\in\ho)
\,\ f_{t-c}\supseteq\r{cl}(\,f_t)\,\}.$
\end{minipage}\end{equation}
(Here we define $f_t$ for all real $t$ by the same
formula~\eqref{d.f_t} that we have been using when $t\in[0,1].$
Thus, in~\eqref{d.D_c}, for $t\leq c$
the set $f_{t-c}$ is all of $X.$
The same set~\eqref{d.D_c} could therefore be defined with
$t$ ranging only over
$(c,1],$
which would avoid the need to extend the definition of $f_t,$
at the price of looking a little less straightforward.)

It is not hard to see that $\Delta_c(X,\r{cl})$
will again form a lattice.
(However, though as noted in \S\ref{S2.cl},
the structure of $\Delta(X,\r{cl})$
is determined by that of the lattice $L_\r{cl},$
the structure of $\Delta_c(X,\r{cl})$ is not determined
by $c$ and $L_\r{cl}.$
Indeed, for $c>0$ its dimension is $\r{card}(X),$
which is not determined by those data.)

Prompted by a comment of the referee's, we ask

\begin{question}\label{Q.thick_ids}
For $X=\{0,1,2\}$ and $\r{cl}$ as in
\S\ref{S2.thickenM3}, and $c\in(0,1),$
what can be said about the lattice
identities satisfied by $\Delta_c(X,\r{cl})$?
In particular, do these depend on $c$?

More generally, if $X$ is any finite set, $\r{cl}$ any
closure operator on $X,$ and $c$ any element of $(0,1),$
what can be said of the variety of
lattices generated by $\Delta_c(X,\r{cl})$?
Will it \textup{(}in general, or in some nice
set of cases\textup{)} have the
property that all finitely generated objects are finite?
\end{question}

There are various further generalizations of this construction.

If we are merely given a family $P$ of subsets of a
set $X,$ then a definition which reduces to~\eqref{d.D_c} when $P$ is
the lattice of closed subsets of a closure operator
is to make $\Delta_c(X,P)$ the poset of those functions
$f{:}\ X\to[0,1]$ such that for all $t>c,$
the set $f_{t-c}$ contains some member of $P$ containing the set $f_t.$

Returning, for simplicity, to the case where $P$ is the
set of closed subsets of $X$ under a closure operator $\r{cl},$
suppose we think of the condition on $f$
in~\eqref{d.D_c} as saying that for each
$Y\subseteq X,$ and $z\in\r{cl}(Y),$
we have $f(z)\geq(\min_{y\in Y}f(y))-c.$
Then one may generalize the construction of~\eqref{d.D_c} to one
in which the $c$ in that
condition is allowed to vary with the pair $(Y,z).$

Finally, we can replace the map $t\mapsto t-c$ by other functions;
for instance, $t\mapsto (1-c)\,t.$

Might some instance of the construction $\Delta_c(X,\r{cl})$
provide a counterexample to Question~\ref{Q.loc_dist}?
No; for let $\r{card}(X)=n,$ and consider the $n!$
bijections $e\,{:}\ \{0,\dots,n{-}1\}\to X.$
Each such bijection determines a simplex in the $\!n\!$--cube $[0,1]^X,$
\begin{equation}\begin{minipage}[c]{25pc}\label{d.D_e}
$D_e\ =\ \{\,f\in[0,1]^X\mid f(e(0))\geq\dots\geq f(e(n{-}1))\},$
\end{minipage}\end{equation}
and the union of these simplices is the whole $\!n\!$--cube.
(These are, in fact, the $\!n\!$--simplices comprising
$\Delta(X,\r{cl'})=[0,1]^X,$ where $\r{cl'}$ is the trivial closure
operation on $X.)$
Let us show that
\begin{equation}\begin{minipage}[c]{25pc}\label{d.D_e_cap}
Each $D_e$ intersects $\Delta_c(X,\r{cl})$
in a distributive sublattice of the latter.
\end{minipage}\end{equation}
This will make $\Delta_c(X,\r{cl})$ a finite union of closed
distributive sublattices, from which it is easy to deduce
a positive answer to Question~\ref{Q.loc_dist} for lattices
$\Delta_c(X,\r{cl}).$

To get~\eqref{d.D_e_cap}, it will suffice to show that for each $e,$
and all $f,g\in D_e\cap\Delta_c(X,\r{cl}),$ we have
\begin{equation}\begin{minipage}[c]{25pc}\label{d.vee_in_D_e_cap}
The coordinatewise supremum and infimum of $f$ and $g$
again lie in $\Delta_c(X,\r{cl}).$
\end{minipage}\end{equation}
This will imply
that these coordinatewise operations give the meet and join
of such a pair $f,\ g$ in $\Delta_c(X,\r{cl}),$ as they
clearly do in $D_e,$ hence they will
carry $D_e\cap\Delta_c(X,\r{cl})$ into itself, making it
a distributive sublattice of $\Delta_c(X,\r{cl}).$

The assertion of~\eqref{d.vee_in_D_e_cap} is straightforward
for the infimum.
To prove the statement about the supremum, consider any
$t\in[0,1],$ any subset $Y\subseteq X,$ and any $z\in\r{cl}(Y).$
We must show that
\begin{equation}\begin{minipage}[c]{25pc}\label{d.sup_fg}
For $f,g\in D_e\cap\Delta_c(X,\r{cl}),$
if the coordinatewise supremum
of $f$ and $g$ is $\geq t$ at all points of $Y,$ then it
is $\geq t-c$ at $z.$
\end{minipage}\end{equation}
Let $i$ be the largest member of $\{0,\dots,n{-}1\}$ such that
$e(i)\in Y.$
If $f$ and $g$ satisfy the hypothesis of~\eqref{d.sup_fg},
their coordinatewise supremum has value at least $t$ at $e(i);$
assume without loss of generality that $f(e(i))\geq t.$
Then by choice of $i,$ and the condition $f\in D_e,$ we
see from~\eqref{d.D_e} that
we have $f(y)\geq t$ for all $y\in Y;$ so $f(z)\geq t-c.$
Hence the coordinatewise supremum of $f$ and $g$ also
has value $\geq t-c$ at $z,$ establishing~\eqref{d.sup_fg},
hence~\eqref{d.vee_in_D_e_cap}, hence~\eqref{d.D_e_cap}, hence
a positive answer to Question~\ref{Q.loc_dist} for
lattices $\Delta_c(X,\r{cl}).$

The above argument goes over to the generalization mentioned
above in which $c$ is allowed to vary with $Y$ and $z.$

\subsection{A construction of D.\,Edmondson}\label{S2.DE}
In an earlier draft of this note, I had a simpler version
of Question~\ref{Q.loc_dist}, asking whether every
topological lattice whose underlying space is a simplicial complex must
be a finite union of closed distributive sublattices.
But it turned out that the example in~\cite{DE} of a nonmodular
lattice structure on the $\!3\!$--cell has properties that
make the existence of such a decomposition very unlikely.
Below, we recall and generalize that construction,
and note why it probably does not admit such a decomposition.

The example of~\cite{DE}, after the change
of coordinates $(x,y,z)\mapsto (x,y,z+xy),$ which turns the
ordering described there into the ordering by componentwise
comparison, is
\begin{equation}\begin{minipage}[c]{25pc}\label{d.DE}
$E\ =\ \{(x,y,z)\in[0,1]^3\mid xy\leq z\leq x\}.$
\end{minipage}\end{equation}

If we write $L=[0,1]^2,$ $M=[0,1],$ each regarded as a
lattice under componentwise operations, and define two set-maps
$L\to M$ by $a(x,y)=xy,$ $b(x,y)=x,$ then the construction
of~\eqref{d.DE} can be abstracted as follows.
(Note that the pair of coordinates $(x,y)$ in~\eqref{d.DE}
and in the above definitions of $a$ and $b$
become the $x$ of~\eqref{d.DEgen}, while the coordinate $z$
of~\eqref{d.DE} becomes the $y$ of~\eqref{d.DEgen}.)

\begin{proposition}[after D.\,Edmondson~\cite{DE}]\label{P.DE}
Let $L$ and $M$ be topological
lattices, and $a, b\,{:}\  L\to M$ continuous isotone
maps \textup{(}not assumed to be lattice homomorphisms\textup{)} such
that $a(x)\leq b(x)$ for all $x\in L.$
Then the set
\begin{equation}\begin{minipage}[c]{25pc}\label{d.DEgen}
$E(L,M;a,b)\ =\ \{(x,y)\in L\times M \mid a(x) \leq y \leq b(x)\},$
\end{minipage}\end{equation}
partially ordered by componentwise comparison, forms
a topological lattice, with operations
\begin{equation}\begin{minipage}[c]{25pc}\label{d.DEops}
$(x,y)\vee(x',y')
\ =\ (x\vee x',\ y\vee y'\vee a(x{\vee}x')),$\\[.5em]
$(x,y)\wedge(x',y')
\ =\ (x\wedge x',\ y\wedge y'\wedge b(x{\wedge}\,x')).$
\end{minipage}\end{equation}
\end{proposition}

\begin{proof}
That the expressions in~\eqref{d.DEops} give a least upper bound and
a greatest lower bound for $(x,y)$ and $(x',y')$ is clear assuming
that the elements described indeed lie in~\eqref{d.DEgen}.
In each case, one of the two inequalities in
the condition for membership in~\eqref{d.DEgen} is
handled by the last joinand or meetand in the
formula; we must in each case verify the
other inequality; namely, for the join we must show that
$y\vee y'\vee a(x{\vee}x')\leq b(x{\vee}x'),$ and for the meet,
that $y\wedge y'\wedge b(x{\wedge}\,x')\geq a(x{\wedge}\,x').$

To get the former inequality, it suffices to verify that each of the
joinands on the left is less than or equal to the term on the right.
This is true of the first joinand, $y,$
because $y\leq b(x)\leq b(x{\vee}x').$
The second joinand, $y',$ is handled in the same way, while the
inequality for $a(x{\vee}x')$ follows from our hypothesis $a\leq b.$
The corresponding result for meets holds by the dual considerations.
Finally, the operations~\eqref{d.DEops} are continuous because
the operations of $L$ and $M,$ and the maps $a$ and $b$
used in~\eqref{d.DEops}, are assumed continuous.
\end{proof}

The condition that the map $a$ of the above construction be
isotone shows that $a(x)\vee a(x')\leq a(x\vee x').$
The first statement of the next result shows that when
that inequality is strict, so that $a$ fails to be a
$\!\vee\!$--semilattice homomorphism, we get a copy
of $N_5$ in our lattice, unless the map $b$ constrains things
too tightly (see~\eqref{d.DE_nonhm}).

\begin{lemma}\label{L.DE_N5}
In the situation of Proposition~\ref{P.DE}, if $x,x'\in L$
are such that
\begin{equation}\begin{minipage}[c]{25pc}\label{d.DE_nonhm}
$a(x)\vee a(x')\ <\ a(x{\vee}x')\wedge b(x),$
\end{minipage}\end{equation}
then the three elements
\begin{equation}\begin{minipage}[c]{25pc}\label{d.DE_N5_gens}
$(x,\,a(x)\vee a(x'))\ <\ (x,\,a(x{\vee}x')\wedge b(x))$
\quad and \quad $(x',\,a(x'))$
\end{minipage}\end{equation}
of $E(L,M;a,b)$ generate a sublattice isomorphic to $N_5.$
\vspace{.2em}

Likewise, if
\begin{equation}\begin{minipage}[c]{25pc}\label{d.DE_nonhm'}
$b(x)\wedge b(x')\ >\ b(x{\wedge}\,x')\vee a(x),$
\end{minipage}\end{equation}
then the three elements
\begin{equation}\begin{minipage}[c]{25pc}\label{d.DE_N5_gens'}
$(x,\,b(x)\wedge b(x'))\ >\ (x,\,b(x{\wedge}\,x')\vee a(x))$
\quad and \quad $(x',\,b(x'))$
\end{minipage}\end{equation}
generate such a sublattice.
\end{lemma}

\begin{proof}
We shall prove the first assertion.
We begin by noting that the first element of~\eqref{d.DE_N5_gens}
indeed satisfies the upper bound on its second
component required by~\eqref{d.DEgen}:
this can be seen from~\eqref{d.DE_nonhm}.
The remaining conditions for our three elements to lie
in $E(L,M;a,b)$ are immediate.

The inequalities $a(x')\leq a(x)\vee a(x')<a(x{\vee}x')\wedge b(x)$
(the last by~\eqref{d.DE_nonhm})
show that on taking the componentwise meet of either of the
first two elements of~\eqref{d.DE_N5_gens} with the third element, we
get $(x\wedge x',\,a(x')),$ hence the corresponding meets
in $E(L,M;a,b)$ are both $(x\wedge x',\,a(x')\wedge b(x{\wedge}x')).$
When we take the corresponding joins, both
likewise give $(x\vee x',\,a(x{\vee}x'))$ (cf.~\eqref{d.DEops}).
So the lattice they generate indeed has the form
\raisebox{2.5pt}[12pt][7pt]{ % N_5
\begin{picture}(15,20)
\dotline{6,-8}{0,-3.5}
\dotline{0,-3.5}{0,3.5}
\dotline{0,3.5}{6,8}
\dotline{6,8}{12,0}
\dotline{12,0}{6,-8}
\end{picture}}.

The second assertion of the lemma holds by the dual observations.
\end{proof}

Turning back to~\eqref{d.DE} and the sentence following it,
which concern the lattice $E$ based on the particular maps
$a,\,b\,{:}\ [0,1]^2\to[0,1]$ given by $a(x,y)=xy,$ $b(x,y)=x,$
observe that for any real numbers $0<x<x'<1$ and $0<y'<y<1$
(note the opposite orders of primed and unprimed
terms in the two inequalities) we have
$a(x,y)\vee a(x',y')=xy\vee x'y',$
while $a((x,y)\vee (x',y'))=a(x',y)=x'y,$ which is strictly larger
than the maximum of $xy$ and $x'y'.$
Moreover, if $(x,y)$ is close enough to $(x',y'),$ then
this strict inequality continues to hold on taking the meet with $b(x)$
as in~\eqref{d.DE_nonhm}.
Thus, by Lemma~\ref{L.DE_N5}
there are little copies of $N_5$ resting all over
the lower surface of the lattice $E.$
This is very different from the behavior of lattices $\Delta(L)$
and the variants of these that we have constructed, where any
nondistributive sublattice that occurs
must have points from more than one simplex.
It appears to me unlikely that the present lattice can, like those,
be written as a finite union of closed distributive sublattices,
as the earlier version of Question~\ref{Q.loc_dist} suggested;
though I don't see a proof that this is impossible.

On the other hand, it is easy to see that every cube
\begin{equation}\begin{minipage}[c]{25pc}\label{d.cube}
$[x,\,x+\varepsilon]\times[y,\,y+\varepsilon]\times[z,\,z+\varepsilon]$
\end{minipage}\end{equation}
lying wholly within our lattice $E$
forms a distributive sublattice; so points of local distributivity
form the whole interior of $E;$ so $E$ does
not give a counterexample to the present form
of Question~\ref{Q.loc_dist}.
(Incidentally, the need for our comment in the preceding paragraph,
that for $(x,y)$ close enough to $(x',y')$ the operation
$-\wedge b(x)$ of~\eqref{d.DE_nonhm} did not matter, can be eliminated
if we take $b(x,y)=1$ in place of $b(x,y)=x$ in~\eqref{d.DE}.
I used the latter only for conformity with \cite{DE}.)

The referee has pointed out that the topological lattice
$\Delta(N_5)$ is itself a case of the construction $E(L,M;a,b).$
Indeed, by the last sentence of \S\ref{S2.cl}, we can identify
$\Delta(N_5)$ with $\Delta(X,\r{cl})$ where $X$ is the set
of join-irreducible elements of $N_5,$
and $\r{cl}$ an appropriate closure operator on $X.$
If we label the join-irreducible elements as follows:
\raisebox{2.5pt}[12pt][7pt]{ % N_5
\begin{picture}(25,20.5)
\dotline{13,-8}{7,-3.5}
\dotline{7,-3.5}{7,3.5}
\dotline{7,3.5}{13,8}
\dotline{13,8}{19,0}
\dotline{19,0}{13,-8}
\put(-1.5,-6){${\scriptstyle 1}$}
\put(19.5,-2.5){${\scriptstyle 2}$}
\put(-1.5,1.5){${\scriptstyle 3}$}
\end{picture}},
then we find that $\Delta(X,\r{cl})$ can be described as
$\{(x_1,x_2,x_3)\in[0,1]^3\mid x_1\wedge x_2\leq x_3\leq x_1\},$
which takes the form $E(L,M;a,b)$ if we let $L=[0,1]^2,$
$M=[0,1],$ $a(x_1,x_2)=x_1\wedge x_2,$ $b(x_1,x_2)=x_1.$
(The two maximal simplices of $\Delta(N_5)$ correspond to
the sets of points with $\!L\!$--coordinate in
$\{(x_1,x_2)\in L\mid x_1\geq x_2\}$ and
$\{(x_1,x_2)\in L\mid x_1\leq x_2\}$ respectively.)

\subsection{Some different sorts of questions}\label{S2.Qnbd}
The above observations on ``tiny'' sublattices of $E(L,M;a,b)$ suggest
the question: To what extent is it true, in a topological lattice, that
``What starts in a small neighborhood stays in a small neighborhood''?
Since a free lattice on three or more generators is infinite,
it seems possible for the sublattice generated by three elements
$p,$ $q,$ $r$ lying in an arbitrarily
small neighborhood of a point $x$ to ``spread'' far from $x.$
An obstruction to this is that this sublattice will lie in the
interval between $p\wedge q\wedge r$ and $p\vee q\vee r,$
each of which will be close to $x\wedge x\wedge x=x\vee x\vee x=x.$
But might the interval between two points that are close to $x$
be very ``wide'', and contain points not close to $x$?

Yes and no.
If $L$ is a compact topological lattice, then an immediate
consequence of~\cite[Lemma~3]{D+S} is that every neighborhood
$U$ of a point $x\in L$ has a subneighborhood $V$ such that
for all $y<z$ in $V,$ all points of the interval
$[y,z]=\{w\in L\mid y\leq w\leq z\}$ lie in $U.$
(This is the case $y<z$ of that lemma, since
when $y<z$ we have $(L\wedge z)\vee y=[y,z].)$
On the other hand, we will see in \S\ref{S2.noncp} that
this can fail in noncompact lattices.

If we look at the sublattice generated by {\em infinitely
many} points in a neighborhood of $x,$ we lose the
trick based on the interval $[p\wedge q\wedge r,\ p\vee q\vee r].$
Here is what little I know about that case.

\begin{lemma}\label{L.nbd}
Let $L$ be a compact topological lattice, and $x$
a point of $L.$
Then the following conditions are equivalent.

\textup{(i)} \ Every neighborhood $U$ of $x$ contains a neighborhood
$V$ of $x$ which is a sublattice of $L.$

\textup{(ii)} \,Every neighborhood $U$ of $x$ contains a neighborhood
$V$ of $x$ which is an interval $[y,z]$ in $L.$
\end{lemma}

\begin{proof}
(ii)$\implies$(i) is immediate, since every interval is a sublattice.

Conversely, assume~(i) holds.
As mentioned above, \cite[Lemma~3]{D+S} gives a
subneighborhood $U'\subseteq U$ of $x$ such that
for all $y<z$ in $U',$ the interval $[y,z]$ lies in $U.$
Because $L$ is compact Hausdorff (recall
Convention~\ref{C.topological}), we can find a {\em closed}
neighborhood $U''\subseteq U'$ of $x;$ and~(i) gives
us a neighborhood $U'''\subseteq U''$
of $x$ which is a sublattice.

The closure of $U'''$ will be a compact sublattice
contained in the closed neighborhood $U'',$
and because it is compact, it will have a least and a greatest
element, $y<z.$
By choice of $U',$ the interval $[y,z]$ is
contained in $U,$ and since it contains $U''',$
it is a neighborhood of $x,$ establishing~(ii).
\end{proof}

One can give further variant statements equivalent to~(i)
and~(ii) above.
In particular, thinking of~(ii) as saying that
every neighborhood $U$ of $x$ contains elements $y\leq z$ such
that $[y,z]$ is a neighborhood of $x,$ we see that it is
equivalent to the
conjunction of the two conditions that $U$ contains
a $y$ such that ${\uparrow}\,y$ is a neighborhood of $x,$ and that it
contains a $z$ such that ${\downarrow}\,z$ is a neighborhood of $x.$

By \cite[Theorem~3]{D+S}, every compact {\em distributive}
topological lattice satisfies condition~(ii) above, hence
both conditions; but for general $L,$ let us ask

\begin{question}\label{Q.near}
If $L$ is a compact Hausdorff topological lattice
\textup{(}perhaps subject to further conditions, such as
being a simplicial complex\textup{),}
must every $x\in L$ satisfy the equivalent conditions
of Lemma~\ref{L.nbd}?
\end{question}

In another direction,
to get a sense of how far things can get from the sort of
examples we have seen, we ask

\begin{question}\label{Q.pathwise}
If $L$ is a connected topological lattice \textup{(}perhaps
assumed compact\textup{)}, must $L$ be pathwise connected?
\end{question}

\subsection{Noncompact lattices}\label{S2.noncp}
We give below a family of examples showing
that much of what the preceding section tells us is true,
or suggests may be true, for compact topological lattices,
fails in the noncompact case.

The idea will be the following.
It is not hard to find lattices $L$ with metrics
$d\,{:}\ L\times L\to\mathbb{R}^{\geq 0}$
in which (say) the interval between a fixed pair of elements $p<q$
contains elements arbitrarily far from $p$ and $q.$
E.g., if $L$ is a lattice consisting of a least and a greatest element,
and an infinite
set of pairwise incomparable elements lying between them, and
$c\,{:}\ L\to\mathbb{N}^{>0}$ is any unbounded function assigning a
positive integer to
each member of $L,$ then letting $d(x,y)=c(x)+c(y)$ when $x\neq y,$
and $0$ otherwise, gives an unbounded metric.
(Think of $c(x)$ as the distance from $x$ to an outside
point $*,$ and imagine the only way to get between distinct
points $x$ and $y$ is via $*.)$
The topology on $L$ defined by this metric is of no interest
to us, since it is discrete; but
the idea will be to start with such a metric on $L,$
and form a lattice of functions $[0,1]\to L,$ metrized
so that appropriate families of functions that differ only on
small subintervals of $[0,1]$ give ``miniaturized'' copies of $L.$
The construction is a bit contrived, but does what we need.
The next lemma describes it; the one that follows
applies it to get the desired examples.

\begin{lemma}\label{L.noncp}
Let $L$ be a lattice,
and $\Gamma(L)$ the quotient of the set of all $\!L\!$--valued
functions on $[0,1]$ that assume only finitely many distinct values,
with each value assumed on a finite union of intervals, by
the equivalence relation that identifies functions which differ
at only finitely many points.
\textup{(}Because of that identification, we don't have to
say whether the abovementioned intervals are open, closed, or
half-open.\textup{)}
Then the pointwise lattice operations determine a lattice
structure on $\Gamma(L)$ under which it
is a subdirect product of copies of $L.$

If, moreover, we are given a real-valued function
$c\,{:}\  L\to[1,\infty),$
such that for all $x,y\in L,$
\begin{equation}\begin{minipage}[c]{25pc}\label{d.c}
$c(x\vee y)\,\leq\,c(x)+c(y)$\quad and\quad
$c(x\wedge y)\,\leq\,c(x)+c(y),$
\end{minipage}\end{equation}
and we let $\Gamma(L,c)$ denote the lattice $\Gamma(L)$ given
with the metric
\begin{equation}\begin{minipage}[c]{25pc}\label{d.d}
$d(\,f,f')\ =\ \int_{\{t\in[0,1]\,\mid\,f(t)\neq f'(t)\}}
\ (c(\,f(t))+c(\,f'(t)))\,dt,$
\end{minipage}\end{equation}
then in the topology induced by that metric, the lattice operations are
continuous, and $\Gamma(L,c)$ is contractible.
\end{lemma}

\begin{proof}
To express $\Gamma(L)$ as a subdirect power of $L,$ let us,
for each $s\in(0,1)$ define $h_s\,{:}\ \Gamma(L)\to L$ to take
each $f\in\Gamma(L)$ to its eventual constant
value as $t$ approaches $s$ from above.
Then we see that each $h_s$ is a homomorphism, and these
homomorphisms separate points.

The verification that for $c$ satisfying~\eqref{d.c},
the formula~\eqref{d.d} defines a metric is straightforward.
(The value is the integral of the $\!(c(x){+}c(y))\!$--metric
of our motivating remarks.)

In proving continuity of the meet and join operations,
it suffices by symmetry to consider joins.
Given $f,\,g\in\Gamma(L,c),$ to prove continuity of $\vee$ at $(\,f,g),$
we must show that
\begin{equation}\begin{minipage}[c]{25pc}\label{d.dvee_need}
For every $\varepsilon>0$ there exists $\delta>0$ such that
for all $f',$ $g'$ satisfying $d(\,f,f')<\delta$ and $d(g,g')<\delta,$
we have $d(\,f\vee g,\,f'\vee g')<\varepsilon.$
\end{minipage}\end{equation}

Now (abbreviating $c(\,f(t))$ to $c(\,f),$ etc., for simplicity),
we see from~\eqref{d.d} and~\eqref{d.c} that
\begin{equation}\begin{minipage}[c]{25pc}\label{d.dvee}
$d(\,f\vee g,\,f'\vee g')\ \leq
\ \int_{\{t\,\mid\,f\neq f'\ {\mbox\scriptsize\r{or}}\ g\neq g'\}}
\ (c(\,f)+c(g)+c(\,f')+c(g'))\,dt.$
\end{minipage}\end{equation}
So to prove~\eqref{d.dvee_need}, it suffices to show that by
taking the $\delta$ of that
statement small enough, we can make the integrals of $c(\,f)+c(\,f')$
and of $c(g)+c(g')$ over the range shown in~\eqref{d.dvee}
each less than $\varepsilon/2.$
By symmetry, it suffices to prove this for $c(\,f)+c(\,f').$
In doing so, let us regard the range of integration as the union of the
set where $f\neq f',$ and the set where $f=f'$ and $g\neq g',$
and show that we can make the integral of $c(\,f)+c(\,f')$
over each of these sets less than $\varepsilon/4.$
The integral over the former set is precisely $d(\,f,f'),$ so we
can make it $<\varepsilon/4$ by taking $\delta\leq\varepsilon/4.$
On the other hand, on the latter set, since $f=f',$
our integrand is $2\,c(\,f).$
Since $f$ is fixed and assumes only finitely many
values, the function $c(\,f)$ has an upper bound $m.$
If we take $\delta\leq\varepsilon/(4m),$ then the condition
$d(g,g')<\varepsilon/(4m),$ implies, in view
of the definition~\eqref{d.d} and
the assumption that $c$ is everywhere $\geq 1,$ that
the total length of the set where $g\neq g'$ must be
$<\varepsilon/(8m).$
Hence the integral of $2\,c(\,f)\leq 2m$ over a subset of
that set will be $\leq\varepsilon/4,$ completing the
proof of continuity.

Finally, to see contractibility, let us, for $f,g\in\Gamma(L,c)$ and
$s\in[0,1],$ define $H_s(\,f,g)\in\Gamma(L,c)$ to agree
with $g$ on $[0,s)$ and with $f$ on $(s,1].$
(Since members of $\Gamma(L,c)$ are defined only up to the relation
of agreeing at all but finitely many points, we don't have to
specify the value at $s$ itself.)
I claim that $H_s(\,f,g)$ is jointly continuous in $s,$ $f,$ and $g.$
Given $f,g,f',g'\in\Gamma(L,c)$ and $s,s'\in[0,1],$ the
triangle inequality in $\Gamma(L,c)$ gives
$d(H_s(\,f,g),\,H_{s'}(\,f',g'))\leq
d(H_s(\,f,g),\ H_{s'}(\,f,g))+d(H_{s'}(\,f,g),H_{s'}(\,f',g')).$
If we fix $s,$ $f$ and $g,$ it is not hard to see that we can
make the distance $d(H_s(\,f,g),\ H_{s'}(\,f,g))$
small by making $s-s'$ small; on the other hand,
we see that $d(H_{s'}(\,f,g),H_{s'}(\,f',g'))$
is bounded by $d(\,f,f')+d(g,g'),$ giving the asserted joint continuity.

In particular, holding $g,$ fixed, $H_s(\,f,g)$ is jointly continuous in
$s$ and $f;$ so letting $s$ vary from $0$ to $1,$ we get
a homotopy from $H_0(-,g),$ which is
the identity function of $\Gamma(L,c),$ to $H_1(-,g),$ which is
the constant function with value $g,$ establishing contractibility.
\end{proof}

We can now give the promised counterexamples.
In the next lemma, note that the condition of the first
sentence, which is assumed in~(i) and~(ii),
holds in particular if $L$ is
an infinite chain with a least and a greatest element,
regarded as a lattice, and $c$
any unbounded function
$L\to[1,\infty).$
On the other hand, the hypothesis of the sentence
containing~(iii) holds for free lattices
on $n\geq 3$ generators, and also for free modular lattices on
$n\geq 4$ generators (see \cite[Exercise~6.3:9]{245}
and \cite[Exercise~I.5.11]{GG.Lat} respectively).

\begin{lemma}\label{L.noncp_cegs}
Suppose that $L$ and $c$ are as in Lemma~\ref{L.noncp},
that $L$ has a least element $0$ and a greatest element $1,$
and that $c$ is unbounded.
Then

\textup{(i)} \  Every nonempty open set $V\subseteq\Gamma(L,c)$
contains elements $p<q$ such that
the interval $[p,q]\subseteq\Gamma(L,c)$ is unbounded \textup{(}i.e.,
has infinite radius under the metric~\eqref{d.d}\textup{)}.

\textup{(ii)}  \,Every nonempty open set $V\subseteq\Gamma(L,c)$
generates a sublattice which is unbounded.
In fact, both the $\!\vee\!$--subsemilattice
and the $\!\wedge\!$--subsemilattice generated by $V$ are
unbounded.

If, rather, $L$ is any infinite lattice generated
by $n<\infty$ elements, and we define $c$ to take each
$x\in L$ to the least length of an expression for $x$ in terms
of those generators \textup{(}where we understand the {\em length}
of an expression to mean the total number of occurrences in that
expression of the symbols for the generators\textup{)}, then

\textup{(iii)} Every nonempty open subset of $V\subseteq\Gamma(L,c)$
contains
$n$ elements which generate an unbounded sublattice of $\Gamma(L,c).$
\end{lemma}

\begin{proof}
To prove~(i),
let us show that for all $r\in\Gamma(L,c)$ and $\varepsilon>0,$ we
can find the desired $p$ and $q$ within distance $\varepsilon$ of $r.$
Suppose we take a small subinterval $I\subseteq [0,1],$ and
let $p$ agree with $r$ {\em except} on $I,$ where it has the value $0,$
and let $q$ likewise agree with $r$ except on $I,$ where
it has the value $1.$
Then if $I$ is small enough (for instance, of length
$<\varepsilon/(2m),$ where $m$ is the maximum of $c(0),$ $c(1),$ and the
images under $c$ of the finitely many values
of the function $r),$ then $p$ and $q$ will
each be within distance $\varepsilon$ of $r.$
Let us now consider any $s\in[p,q]$ which likewise agrees with $r$
except on $I,$ but has for value there
an arbitrary element $x\in L.$
If we take $x$ such that $c(x)$ is sufficiently
large, then $s$ will be arbitrarily far from $r;$ so
$[p,q]$ is indeed unbounded.

To prove the $\!\vee\!$--semilattice case of~(ii), we shall show
that for any $r\in\Gamma(L,c),$ $V$ contains a
finite family of points $p_i$ such that $\bigvee p_i$
is ``far'' from $r.$
The construction is similar to the preceding, so I shall
abbreviate the details.
We again start by replacing the values of $r$ on a small subinterval
$I\subseteq [0,1]$ with $0\in L,$ getting an element $p\in\Gamma(L,c).$
We then take $x\in L$ such that $c(x)$ has large value, and
let $N$ be the integer $\lceil c(x)\rceil.$
We subdivide $I$ into $N^2$ equal small
subintervals $I_i$ $(0\leq i< N^2),$
and let each $p_i$ be formed from $p$
by changing the value on $I_i$ from $0$ to $x.$
Since the length of each $I_i$ is $1/N^2$ times
that of $I,$ the distance $d(p,p_i)$ is $\leq (N+c(0))/N^2$
times the length of $I,$ hence is small if $c(x),$ and hence $N,$
is large enough; so the $p_i$ can all be made close to $p,$
which, if $I$ has been taken small, is close to $r.$
But if we form the join of these $N^2$ elements in $\Gamma(L,c),$
this has value $x$ on {\em all} of $I;$ so for fixed $I$
and large enough $c(x),$ this join is arbitrarily
far from $p,$ and hence from $r.$
For the $\!\wedge\!$--semilattice case of~(ii) we use the dual
construction.

Turning to~(iii), it is easy to see that the function $c$ defined
in the sentence containing that
statement satisfies~\eqref{d.c} and is unbounded.
To get an unbounded sublattice of $\Gamma(L,c)$ generated by $n$
elements near an element $r\in\Gamma(L,c),$ again
choose a small interval $I\subseteq [0,1],$ and this time construct
$n$ elements by replacing the values of $r$ on $I$ by each
of our $n$ generators of $L.$
It is easy to see that for $I$ small enough,
those $n$ elements all lie in $V,$ but
that no matter how small $I$ is, the unboundedness of $c$
implies that the sublattice they
generate has elements arbitrarily far from $r.$
\end{proof}

By choosing $L$ as indicated immediately before the
above lemma, we can get distributive $\Gamma(L,c)$ satisfying
conditions~(i) and~(ii), and modular $\Gamma(L,c)$ satisfying~(iii).

The construction $\Gamma(L,c)$ really uses
nothing specific to lattices, and might have other
applications in universal algebra, if it is not already known.

\section{Further notes on the order-complex construction}\label{S.further}
\subsection{Why topologists haven't looked at \texorpdfstring{$\Delta(L)\!$}{Delta(L)}}\label{S2.L-01}
In the literature I am aware of, e.g.,~\cite{MW},
the order complex construction
$\Delta(\ )$ is regularly applied, not to finite lattices $L,$ but to
their subposets $L-\{0,1\}.$
This is because $\Delta(L-\{0,1\})$ can
have nontrivial homotopy and homology, while $\Delta(L)$ cannot,
nor can $\Delta(L-\{0\})$ or $\Delta(L-\{1\}).$
Indeed,

\begin{lemma}\label{L.contractible}
If a finite partially ordered set $P$ has a least element,
a greatest element, or more generally, an element $z$ that is
comparable to all elements of $P,$ then $\Delta(P)$ is contractible.
\end{lemma}

\begin{proof}
Given $z$ comparable to all elements of $P,$ let us
define for each $s\in[0,1]$ a map $H_s\,{:}\ \Delta(P)\to\Delta(P);$
namely, for $f\in\Delta(P),$ let
$H_s(\,f) = (1-s)\cdot f + s\cdot({\downarrow}\,z).$
To see that $H_s(\,f)$ lies in $\Delta(P),$
note that $f$ is a convex linear
combination of characteristic functions of principal ideals
corresponding to
a chain of elements of $P,$ and since $z$ is comparable with
all elements of $L,$ throwing it in still leaves a chain.
Since $H_0$ is the identity map of $\Delta(P)$
and $H_1$ is a constant map (i.e., $H_1(\,f)={\downarrow}\,z$ is
independent of $f),$ $\Delta(P)$ is contractible.
\end{proof}

Examining the above proof, one sees that the simplicial
complex $\Delta(P)$ is a cone on $\Delta(P-\{z\}).$
In particular, if $L$ is a finite lattice with more than one element,
then $\Delta(L)$ is a {\em cone on a cone} on $\Delta(L-\{0,1\}).$
For instance, if $L=M_3,$ then $L-\{0,1\}$ is a $\!3\!$--element
antichain, hence $\Delta(L-\{0,1\})$ is a $\!3\!$--point discrete space.
The simplicial complex $\Delta(L-\{0\})$ or $\Delta(L-\{1\})$
is a cone on that set, and so looks like
\raisebox{2pt}[12pt][8pt]{ % oY
\begin{picture}(20,20)
\qbezier(6,7)(6,7)(6,0)
\qbezier(0,-4)(0,-4)(6,0)
\qbezier(12,-4)(12,-4)(6,0)
\end{picture}}
(with the center point as greatest or least element, respectively);
while $\Delta(L),$ a cone on the latter space, has, as we have
seen, the form of three triangles connected along a common edge.

Another triviality result, specific to lattices, is

\begin{lemma}\label{L.retract}
If $L$ is a finite lattice, then $\Delta(L)$ is a retract
of $[0,1]^L$ by a piecewise linear order-preserving map.
More generally, if $\r{cl}$ a closure operator on a finite
set $X,$ then $\Delta(X,\r{cl})$ is a retract
of $[0,1]^X$ by such a map.
\end{lemma}

\begin{proof}[Sketch of proof]
To get a retraction $r$ in the more general situation,
map each $f\in[0,1]^X$ to the function $r(\,f)$ such that for
$t\in\ho,$ $r(\,f)_t=\r{cl}(\,f_t);$ equivalently, such that for
$x\in L,$ $r(\,f)(x)$ is the largest $t$ such that $x\in\r{cl}(\,f_t).$
\end{proof}

For another such result, see \cite[Theorem~6.2]{WT77}.

Neither of the above lemmas holds for the class of lattices of the
form $\Delta(L,S)$ described in Proposition~\ref{P.pairs}.
For example, if $L=\{0,1\}$ and $S=\{(1,1)\},$ then
$\Delta(L,S)$ is a $\!2\!$--point disconnected lattice.
However, if we exclude pairs of the form $(x,x)$ from $S,$ then
the result of Lemma~\ref{L.contractible} does go over, though
it takes more work to prove:

\begin{lemma}\label{L.*D(L,S)_cnctd}
Let $L$ be a nontrivial finite lattice, and $S\subseteq L\times L$
a family of pairs $(x,y)$ all of which satisfy $x<y$
\textup{(}cf.\ the weaker condition $x\leq y$
of Proposition~\ref{P.pairs}\textup{)}, and have $y$ join-prime
\textup{(}as in Proposition~\ref{P.pairs}\textup{(iii))}.

Then for every minimal nontrivial interval $[p,q]$ of $L$
\textup{(}i.e., for every pair of elements $p<q$ in $L$
with no element $r$ strictly between $p$ and $q),$
$\Delta(L,S)$ contains the edge from ${\downarrow}\,p$
to ${\downarrow}\,{q}$ in $\Delta(L).$

Hence $\Delta(L,S)$ is pathwise connected.
Hence it is, in fact, contractible.
\end{lemma}

\begin{proof}
Let $[p,q]$ be a minimal nontrivial interval in $L.$
The points of the edge from ${\downarrow}\,{p}$
to ${\downarrow}\,{q}$ in $\Delta(L)$ have the form
\begin{equation}\begin{minipage}[c]{25pc}\label{d.p_to_q}
$f\ =\ (1-t)\,({\downarrow}\,p)+t\,({\downarrow}\,{q})
\quad (t\in[0,1]).$
\end{minipage}\end{equation}
By the definition of $\Delta(L,S),$ if $\Delta(L,S)$
does not contain the above edge, there must be $(x,y)\in S$
such that some element~\eqref{d.p_to_q}
carries both  $x$ and $y$ to values in $(0,1);$
which means that both $x$ and $y$ lie
in $({\downarrow}\,q)-({\downarrow}\,p).$
This implies $p < p\vee x \leq q;$ so as
the interval $[p,q]$ is minimal, we must have $p\vee x = q.$
But the relation $y\leq q=p\vee x$ contradicts the
join-primeness of $y.$
This contradiction gives the first assertion of the lemma.

That assertion says, in effect, that $\Delta(L,S)$ contains
the Hasse diagram of the lattice $L.$
But that diagram connects all the vertices of $\Delta(L),$
equivalently, of $\Delta(L,S).$
So $\Delta(L,S)$ is pathwise connected, as claimed.

In particular, there exists a path $h\,{:}\ [0,1]\to\Delta(L,S)$
with $h(0)={\downarrow}\,0$ and $h(1)={\downarrow}\,1.$
Starting with such an $h,$ define
$C\,{:}\ [0,1]\times\Delta(L,S)\to\Delta(L,S)$
by $C(t,f)(z)=h(t)\vee f(z).$
Then $C$ is continuous,
$C(0,-)$ is the identity function of $\Delta(L,S),$
and $C(1,-)$ is the constant function taking all
$f\in\Delta(L,S)$ to ${\downarrow}\,1.$
This proves contractibility.
\end{proof}

\subsection{\texorpdfstring{$\Delta(L)$}{Delta(L)} as a subspace of \texorpdfstring{$L^{[0,1]}$}{L sup [0,1]}}\label{S2.subdirect}
We saw in \S\ref{S2.M} that the
construction $\Delta(L)$ is the case $M=[0,1]$ of a general
construction that is symmetric in two complete lattices $L$ and $M.$
Hence, though we have mainly been regarding $\Delta(L)$ as a subspace
of $[0,1]^L,$ the symmetric description shows that
it can also be viewed as a subspace of $L^{[0,1]}.$
Namely, to each $f\in\Delta(L)\subseteq[0,1]^L$ we can
associate the member of $L^{[0,1]}$ taking each $t\in[0,1]$
to the element $h_t(\,f)\in L$ such that $f_t={\downarrow}\,h_t(\,f).$
(Cf.\ last two sentences of first paragraph of Theorem~\ref{T.main}.)

Now for $L$ a finite lattice, $L^{[0,1]},$ being
an infinite direct product of finite sets, has a natural
totally disconnected compact Hausdorff topology.
This induces a topology on $\Delta(L),$ in general stronger
(having more open sets) than its topology as a simplicial complex.
To see that it is at least as strong as the latter, note
that given any $\varepsilon>0,$ if we choose $0<t_1<\dots<t_n<1$
in $[0,1]$ with successive differences $<\varepsilon,$ then
the coordinates indexed by $t_1,\dots,t_n$ of the
image in $L^{[0,1]}$ of an element $f\in\Delta(L)$
determine to within $\varepsilon$ all
the coordinates of $f$ as a member of $[0,1]^L.$
Thus, in $\Delta(L),$ every neighborhood of $f$
with respect to the $\![0,1]^L\!$--topology
contains a neighborhood of $f$
with respect to the $\!L^{[0,1]}\!$--topology.
These topologies on $\Delta(L)$ are distinct
if $L$ has more than one element, since one
is totally disconnected and the other is contractible.
The $\!L^{[0,1]}\!$--topology on $\Delta(L)$
is not, of course, compact, since
distinct compact Hausdorff topologies on a space cannot be comparable.
It follows that $\Delta(L),$ though
closed in $[0,1]^L,$ is not, in general, closed in~$L^{[0,1]}.$

To see concretely the difference between the
two topologies, note that if one moves
$f\in\Delta(L)$ continuously in the $\![0,1]^L\!$--topology, so that
its value at some $x\in L$ comes up from below to a value
$t\in\ho,$ then when it reaches that value,
the image of $f$ in the $\!t\!$--indexed
coordinate of $L^{[0,1]}$ jumps discontinuously from an element
of $L$ that does not majorize $x$ to one that does.
(However, we do have ``continuity from above'':
if we move $f$ {\em downward} continuously in $\Delta(L)$
with respect to the $\![0,1]^L\!$--topology, then
its image in $L^{[0,1]}$ at any
moment {\em is} the limit of its earlier values.)

The failure of $\Delta(L)$ to be closed in $L^{[0,1]},$ deduced
above, can be seen directly.
Let $L$ be the two-element lattice $\{0,1\}.$
Then $\Delta(L),$ as a subset of $L^{[0,1]}=\{0,1\}^{[0,1]},$
corresponds to the set of characteristic functions of the principal
ideals of $[0,1],$ i.e., the characteristic
functions of the subintervals $[0,t]$ $(0\leq t\leq 1).$
It is easy to see that the closure of this set in the
\mbox{$\!L^{[0,1]}\!$--topology} consists of these functions,
together with the characteristic functions of the nonempty
{\em non-principal} ideals, i.e., the sets
$[0,t)$ $(0<t\leq 1).$

Under both topologies,
the meet and join operations of $\Delta(L)$ are continuous.

In subsequent sections we will, as before, consider
our lattices $\Delta(L)$ to have the $\![0,1]^L\!$--topology
(and similarly for the objects $\Delta(P)$ where
$P$ is a finite partially ordered set or a $\!\vee\!$-
or $\!\wedge\!$--semilattice).

\subsection{Functoriality}\label{S2.funct}
Like the familiar power-set construction on sets, our construction
$\Delta(\ )$ on posets, semilattices and lattices can
be made into both covariant and contravariant functors; though in this
case we shall find that the contravariant constructions are more
limited and less straightforward than the covariant ones.
We begin with the covariant case.

\begin{proposition}\label{P.cov_funct}
The construction $P\mapsto\Delta(P)$ of \S\ref{S.constr} can be
made a covariant functor from the category of finite posets and isotone
maps to the category of partially ordered simplicial complexes
and isotone simplicial maps.
Namely, given an isotone map $h\,{:}\ P\to Q,$ we define the
induced map $\Delta(h)\,{:}\ \Delta(P)\to\Delta(Q)$ to carry an
element $f$ which is a convex linear combination of characteristic
functions ${\downarrow}\,x_1,\dots,{\downarrow}\,x_n$ $(x_i\in P)$
to the convex linear combination
of ${\downarrow}\,h(x_1),\dots,{\downarrow}\,h(x_n)$ with the
same coefficients.
Equivalently, for each $y\in Q,$ $(\Delta(h)(\,f))(y)$ is the maximum
of $f(x)$ over all $x\in P$ with $h(x)\geq y.$

When applied to posets $P$ in which every pair of elements
has a join \textup{(}respectively, a meet\textup{)},
and morphisms $h$ that respect joins \textup{(}meets\textup{)},
the functor gives morphisms $\Delta(h)$ that likewise respect joins
\textup{(}meets\textup{)}.
Thus, it also gives covariant functors from the categories of finite
$\!\vee\!$--semilattices,  $\!\wedge\!$--semilattices, and
lattices to those of topological
$\!\vee\!$--semilattices, $\!\wedge\!$--semilattices and lattices.
\end{proposition}

\begin{proof}[Sketch of proof]
It is straightforward to see that the construction described in the
first paragraph gives a functor, and has the equivalent description
noted in the last sentence of that paragraph.

In view of the descriptions of meets and joins in $\Delta(P)$
given in Theorem~\ref{T.main}, the proof of the remaining assertions
comes down to showing that if $P$ has joins,
respectively meets, and $Q$ does likewise, and these are
respected by $h,$ then for all $f,\,g\in\Delta(P)$ and
$t\in\ho,$ we have
$(\Delta(h)(\,f\vee g))_t=(\Delta(h)(\,f))_t\vee (\Delta(h)(g))_t,$
respectively
$(\Delta(h)(\,f\wedge g))_t =
(\Delta(h)(\,f))_t\wedge (\Delta(h)(g))_t,$ as principal ideals of $Q.$
Now in both $P$ and $Q,$ joins (respectively,
meets) of principal ideals correspond to the joins (respectively,
meets) of their generating elements.
Hence if $h$ respects joins and/or meets of elements of
$P$ and $Q,$ the map $\Delta(h),$ as defined
in the statement of the proposition, will respect these operations
on elements of $\Delta(P)$ and $\Delta(Q).$
\end{proof}

We now give the more limited result holding
for $\Delta(\ )$ as a contravariant functor.
In connection with the statement below, note that in a {\em finite}
$\!\vee\!$--semilattice $P$ with $0,$ every subset has not only
a least upper bound, but also a greatest lower bound;
and dually for finite $\!\wedge\!$--semilattices with $1.$
Hence the difference between such semilattices and finite
lattices merely concerns which maps among them are considered morphisms.
In all these cases, the existence of
both least upper bounds and greatest lower bounds necessarily
carries over to the ordered simplicial complex $\Delta(P),$
since the construction of that object is the same whether we
regard $P$ as a poset, semilattice, or lattice.

\begin{proposition}\label{P.contrav_funct}
On the category whose objects are finite partially ordered
sets, and whose morphisms are set-maps $h$ satisfying
\begin{equation}\begin{minipage}[c]{25pc}\label{d.prin-prin}
The inverse image under $h$ of every principal ideal is a
principal ideal,
\end{minipage}\end{equation}
\textup{(}a subset of the isotone maps\textup{)},
the construction $P\mapsto\Delta(P)$ of \S\ref{S.constr} can be
made a {\em contravariant} functor $\Delta'$
to the category of partially ordered
simplicial complexes and isotone simplicial maps.
Namely, for $h\,{:}\ P\to Q$ satisfying~\eqref{d.prin-prin}, we let
$\Delta'(h)\,{:}\ \Delta(Q)\to\Delta(P)$ carry $f\in\Delta(Q)$ to
$f\circ h\in\Delta(P).$

In particular, $\Delta'$ yields a contravariant functor
from the category of finite $\!\vee\!$--semilattices
with least element, $0,$ and morphisms $h$ preserving joins
and least elements, to that of partially ordered simplicial complexes.
However, the resulting morphisms $\Delta'(h)$ need {\em not}
respect the join operations inherited by the complexes $\Delta(P)$
\textup{(}as described in Theorem~\ref{T.main} or
Proposition~\ref{P.cov_funct}\textup{)}, even when applied to
$\!0\!$-- and $\!1\!$--respecting lattice homomorphisms between
finite distributive lattices;
but {\em will} respect the {\em meet} operations and
greatest elements $1$ of these complexes,
even if the given morphisms of $\!\vee\!$--semilattices do not.
\end{proposition}

\begin{proof}
To see that a map $h\,{:}\ P\to Q$ satisfying~\eqref{d.prin-prin}
must be isotone, let $x\leq y$ in $P.$
Then $h^{-1}({\downarrow}\,h(y))$ will be a principal ideal
of $P$ containing $y,$ hence containing $x\leq y,$ which says that
$h(x)\in{\downarrow}\,h(y),$ i.e., $h(x)\leq h(y),$ as required.

Given such an $h\,{:}\ P\to Q,$ condition~\eqref{d.prin-prin}
shows that if $f{:}\ Q\to[0,1]$ belongs to $\Delta(Q),$ then the element
\begin{equation}\begin{minipage}[c]{25pc}\label{d.*D'hf}
$\Delta'(h)(\,f)\ =\ f\circ h$
\end{minipage}\end{equation}
belongs to $\Delta(P).$
It is straightforward to verify that this makes $\Delta'(h)$ an
isotone simplicial map, and that the construction $\Delta'$ so defined
respects (contravariantly)
composition of morphisms and identity morphisms.

To see that when $P$ and $Q$ are $\!\vee\!$--semilattices
with $0,$ every map $h\,{:}\ P\to Q$ respecting those operations
satisfies~\eqref{d.prin-prin}, note that the inverse image
under $h$ of any principal ideal of $Q$ is nonempty
(because it contains $0_P),$
is a downset, and is closed under $\vee;$ hence, since
$P$ is finite, it will be a principal ideal.

For an example showing that even for $h$ a homomorphism of finite
distributive lattices which respects greatest and least
elements, the map $\Delta'(h)$ need not respect joins,
let $Q$ be the $\!8\!$--element lattice $\{0,1\}^3,$ $P$ its sublattice
consisting of the four elements $(i,j,0)$ and the one
element $(1,1,1),$ and $h$ the inclusion of $P$ in $Q.$
Then the characteristic functions of ${\downarrow}(1,0,1)$ and
${\downarrow}(0,1,1)$ in $\Delta(Q)$ have for inverse images in
$\Delta(P)$ the characteristic functions of ${\downarrow}(1,0,0)$ and
${\downarrow}(0,1,0),$ whose join is ${\downarrow}(1,1,0).$
However, the join of the given characteristic functions in $\Delta(Q)$
is ${\downarrow}(1,1,1),$ whose inverse image in $\Delta(P)$ is
${\downarrow}(1,1,1).$

One sees, however, that the maps $\Delta'(h)$ respect meets,
including the empty meet, by noting
that $(\,f\circ h)_t=h^{-1}(\,f_t),$ and combining this with the facts
that taking inverse images under a set map respects intersections
of subsets (including the empty intersection, i.e.,
the greatest element), and that intersection as subsets gives
the meet operation on principal ideals in lattices.
\end{proof}

\subsection{Representability}\label{S2.rep}
The covariant version of our
order complex construction,
and to a more limited extent, the contravariant version,
are {\em representable}
\cite[Definition in {\S}III.2]{CW} \cite[Definition~8.2.3]{245}.

The superficially more obvious case is the contravariant one:
For $P$ a finite partially ordered set, $\Delta(P)$
is a collection of antitone (i.e., order-reversing) maps from $P$ to
the interval $[0,1],$ so its points should be the morphisms
$P\to[0,1]^{\r{op}}$ in an appropriate category.
If we denote by ${\bf Poset}_{\r{prin{\hy}id}}$
the overcategory
of the category of Proposition~\ref{P.contrav_funct},
whose objects are now {\em all} partially
ordered sets, but whose morphisms are still the maps $h$
which satisfy~\eqref{d.prin-prin}, then this works --
except where the anomalous behavior of $0\in [0,1]^\r{op},$
discussed in~\S\ref{S2.01}, interferes.
We find that
\begin{equation}\begin{minipage}[c]{25pc}\label{d.*D'=}
${\bf Poset}_{\r{prin{\hy}id}}(P,\,[0,1]^\r{op})\ =
\left\{ \begin{array}{cl}
\Delta(P) & \mbox{if $P$ has a greatest element,}\\[.2em]
\emptyset & \mbox{otherwise}.
\end{array}\right.$
\end{minipage}\end{equation}

We shall not look further at this contravariant functor.
(The reader can verify that replacing $[0,1]^\r{op}$ with $\ho^\r{op}$
does not cure the ``$\emptyset$ otherwise'' feature of~\eqref{d.*D'=}.)

For the covariant functor $\Delta,$ things work out more nicely.
We saw in \S\ref{S.constr} that each $f\in\Delta(P)$
is determined by the chain of principal ideals $f_t$
as $t$ ranges over $\ho.$
Each such principal ideal $f_t$
has the form ${\downarrow}\,h_t(\,f)$ for some
$h_t(\,f)\in P$ (\S\ref{S2.Delta}), so $f$ is determined by
this chain of elements, which we can regard as
a morphism $h(\,f)\,{:}\ \ho^{\r{op}}\to P.$

In what category?
The above category ${\bf Poset}_{\r{prin{\hy}id}}$ would do;
but in fact, the morphisms $\ho^{\r{op}}\to P$
in that category are the same as the
morphisms $\ho^{\r{op}}\to P$ in the much larger category
${\bf Poset}_{\r{chain{\hy}sup}}$ of isotone maps that
respect least upper bounds of chains (when these exist).
So for finite partially ordered sets $P,$ we can make
the identification
\begin{equation}\begin{minipage}[c]{25pc}\label{d.*D=}
${\bf Poset}_{\r{chain{\hy}sup}}(\kern.08em\ho^\r{op},\,P)
\ =\ \Delta(P).$
\end{minipage}\end{equation}
Note that isotone maps among finite partially
ordered sets trivially respect least upper bounds of chains,
explaining why in Proposition~\ref{P.cov_funct} we were able
to define the covariant version of $\Delta$ on all such maps.
It is straightforward to verify that the behavior of
${\bf Poset}_{\r{chain{\hy}sup}}(\ho^\r{op},\,-)$
on morphisms agrees with the behavior described there; in
particular, it gives
isotone maps that are continuous, and in fact are simplicial maps.

Covariant representable set-valued functors
are known to preserve {\em limits}
in the category-theoretic sense (products of objects, equalizers
of pairs of morphisms, and constructions that can be obtained
from these) \cite[Theorem~V.4.1]{CW}, \cite[Theorem~8.8.7]{245}.
Let us show that $\Delta$ respects limits, not only
as a set-valued functor, but as a functor to partially
ordered topological spaces.
Since we are only considering $\Delta$ to be
defined on {\em finite} partially
ordered sets, our result will concern limits over finite diagrams.

\begin{theorem}\label{T.*D_resp_lim}
As a covariant functor from finite partially ordered sets to
partially ordered compact Hausdorff spaces, $\Delta$
respects limits over finite diagrams.

Hence the same is true of $\Delta$ as a covariant functor
from finite $\!\vee\!$--semilattices, $\!\wedge\!$--semilattices,
or lattices, to compact Hausdorff spaces with structures of topological
$\!\vee\!$--semilattice, $\!\wedge\!$--semilattice, or lattice.
\end{theorem}

\begin{proof}[Sketch of proof]
Let ${\bf FinPoset}$ be the category of finite partially ordered
sets and isotone maps, which we have noted is a full
subcategory of ${\bf Poset}_{\r{chain{\hy}sup}}.$
Let ${\bf D}$ be a finite category
(one with finitely many objects and finitely many morphisms),
and $F{:}\ {\bf D}\to{\bf FinPoset}$ a functor.
Since both ${\bf FinPoset}$ and the category of partially
ordered compact Hausdorff spaces admit finite limits, we get a
comparison morphism
\begin{equation}\begin{minipage}[c]{25pc}\label{d.C_compar}
$C\,{:}\ \ \Delta(\varprojlim_{\bf D} F)\ \to
\ \varprojlim_{\bf D}(\Delta\circ F)$
\end{minipage}\end{equation}
\cite[Definition~8.8.13 et seq.]{245}.

Because at the set level $\Delta$ is representable
within the category ${\bf Poset}_{\r{chain{\hy}sup}},$
and therefore respects limits, the continuous map~\eqref{d.C_compar}
will be a bijection.
But a continuous bijection between compact Hausdorff spaces
is a homeomorphism.

To show that $C$ is an order isomorphism,
let us again call on the fact that within the larger category
${\bf Poset}_{\r{chain{\hy}sup}},$ the object $\ho^\r{op}$
represents the functor $\Delta$ at the underlying set level.
It follows that if we form the object \mbox{$\ho^\r{op}\times\{0,1\},$}
with the product ordering, then a morphism
from this object to a finite partially ordered set $P$
will correspond to a pair $(\,f,g)$ of members of $\Delta(P)$
with $f\leq g;$ hence $\ho^\r{op}\times\{0,1\}$ represents the functor
giving the graph of the order relation on $\Delta(P).$
Since this functor is representable, it respects
limits; so the orderings on $\Delta(\varprojlim F)$ and
$\varprojlim(\Delta\circ F)$ agree, as claimed.

The final sentence holds because limits of the indicated
algebraic structures have as their underlying partially
ordered sets, respectively
topological spaces, the limits of those sets or spaces.
\end{proof}

The case of the above result where ${\bf D}$ is a $\!2\!$--object
category with only the two identity morphisms
says that for any two finite posets $P$ and $Q$ we have
\begin{equation}\begin{minipage}[c]{25pc}\label{d.P1xP2}
$\Delta(P\times Q)\ \cong\ \Delta(P)\times\Delta(Q).$
\end{minipage}\end{equation}
This isomorphism is classical
\cite[Lemma~II.8.9]{E+S}, \cite[Theorem~3.2]{JW},
but the proofs in those works are quite complicated.
In~\S\ref{S2.PxQ} we will give a ``hands-on'' proof of the
same result, which still seems a bit simpler than those cited.

Can we say anything about the behavior of $\Delta$ on {\em colimits}?
In general, covariant representable functors do not respect colimits;
but in some particular classes of cases they do \cite[\S8.9]{245}
\cite{dirlimfix}.
I claim that our functor
$\Delta$ respects finite coproducts of partially ordered sets.
Note that a coproduct $P\smallcoprod Q$ of posets is
their disjoint union, with each poset
having its original order structure,
and no order-relations holding between
elements of $P$ and elements of $Q.$
Hence, an isotone map from $\ho^\r{op}$ (or, indeed, any
chain) to $P\smallcoprod Q$
must have image entirely in $P$ or entirely in $Q.$
So $\Delta(P\smallcoprod Q)$ is a disjoint union of
$\Delta(P)$ and $\Delta(Q);$ and it is easy to see that it also has the
right structure of ordered topological space to be their coproduct.

But this does not yield corresponding statements for
coproducts of semilattices or lattices.
Colimits of algebraic structures, unlike limits, do
not typically arise from colimits of their underlying sets
(or, in this case, their underlying posets), so the above
argument does not imply that $\Delta$ respects coproducts of
lattices or semilattices; and in fact, it does not.
For instance, if $P$ and $Q$ are each $\!1\!$--element lattices,
then their coproduct as lattices is the $\!4\!$--element lattice
\raisebox{2.5pt}[12pt][7pt]{ % <>
\begin{picture}(15,20)
\dotline{7,8}{0,0}
\dotline{0,0}{7,-8}
\dotline{7,-8}{14,0}
\dotline{14,0}{7,8}
\end{picture}}\,,
and applying $\Delta$ to that, we get a lattice consisting
of two triangles $\!(2\!$--simplices) fitted together to form a diamond.
In contrast, the topological lattices $\Delta(P)$ and
$\Delta(Q)$ are, like $P$ and $Q,$ $\!1\!$--point lattices,
whose coproduct, as for discrete lattices, has just $4$ elements.

\subsection{A discrete analog of \texorpdfstring{$\Delta(P)\cong\Delta'(P)$}{Delta(P)=Delta(Q)}}\label{S2.discrt_contra}

The fact that the geometric construction $\Delta(P)$ can be
regarded as the values of either a covariant or a contravariant functor
has an analog purely in the realm of discrete posets $P,$ in which,
moreover, there is no reason to restrict attention to finite $P.$
I sketch it below.

The idea is to replace the interval $[0,1]$ in the construction $\Delta$
with the $\!2\!$--element object $\{0,1\}.$
Thus, we associate to every poset $P$ the set of
$\!\{0,1\}\!$--valued functions $f{:}\ P\to\{0,1\}$ such that
the set $f^{-1}(1)$ is a principal ideal.
This, of course, is just the set of characteristic functions
of principal ideals of $P,$ and can be identified with $P;$
and indeed, if we make this construction a {\em covariant} functor as in
Proposition~\ref{P.cov_funct}, it is the
identity functor of ${\bf Poset}.$

The {\em contravariant} functor analogous
to $\Delta'$ should be a functor on ${\bf Poset}_{\r{prin{\hy}id}},$
carrying each morphism $h\,{:}\ P\to Q$ in that category
to the map that takes every principal ideal ${\downarrow}\,q$
to the principal downset $h^{-1}({\downarrow}\,q)$ of $P;$
or, looked at as a map on elements, that takes
$q$ to the generator of $h^{-1}({\downarrow}\,q).$
Let us call this functor~$G.$

As in the case of $\Delta',$ the maps $G(h)\,{:}\ Q\to P$
will be isotone.
However, they will not in general
be morphisms of ${\bf Poset}_{\r{prin{\hy}id}}.$
Rather, they turn out to have the dual property
of carrying each principal {\em filter}
${\uparrow}\,p$ to a principal filter of $Q,$
namely ${\uparrow}\,h(p).$
In fact, one finds that $G$ gives a contravariant
equivalence between ${\bf Poset}_{\r{prin{\hy}id}}$ and
${\bf Poset}_{\r{prin{\hy}filt}}.$

Given any two partially ordered sets $P$ and $Q,$ the structure
given by a pair of morphisms related under this equivalence,
\begin{equation}\begin{minipage}[c]{25pc}\label{d.Gal}
$P \xrightarrow[{\bf Poset}_\r{prin{\hy}id}]{h}\ Q$ \qquad and \qquad
$Q \xrightarrow[{\bf Poset}_\r{prin{\hy}filt}]{G(h)}\ P$
\end{minipage}\end{equation}
is what is called a {\em Galois connection}
between $P$ and $Q,$ with $h$ the ``lower adjoint''
and $G(h)$ the ``upper adjoint''~\cite{wiki_Gal}.

The functor $\Delta'$ of \S\ref{S2.funct}
can be described in terms of the above construction,
as the composite $\Delta\circ G.$

If $P$ and $Q$ are complete lattices, it is not hard to show
that the morphisms $P\to Q$ in ${\bf Poset}_{\r{prin{\hy}id}}$ are
precisely the maps $P\to Q$ that respect all joins
(empty, finite, and infinite),
while morphisms in ${\bf Poset}_{\r{prin{\hy}filt}}$ are
those that similarly respect all meets.
Hence the functor $G$ yields a contravariant
equivalence between the category of complete lattices and maps
that respect arbitrary joins -- in other words, the category of
complete $\!\vee\!$--semilattices with least element, with morphisms
respecting all the operations (zeroary, finitary and infinitary)
-- and the category with the same objects, but
morphisms similarly respecting all meets -- the category of
complete $\!\wedge\!$--semilattices with greatest element.
Restricting attention to finite lattices, and composing
with $\Delta,$ this explains the
properties of $\Delta'$ noted in the last paragraph
of Proposition~\ref{P.contrav_funct}.

\subsection{A hands-on construction of \texorpdfstring{$\Delta(P\times Q)$}{Delta(P x Q)}}\label{S2.PxQ}
Returning to the geometric construction $\Delta(-),$
let us give, as promised, a proof of the
case of Theorem~\ref{T.*D_resp_lim} where the limit in question
is a pairwise direct product, which does not rely on general
properties of representable functors, but only on
the description of $\Delta$ as a functor in \S\ref{S2.funct}.

\begin{proposition}[{\cite[Lemma~II.8.9]{E+S}, \cite[Theorem~3.2]{JW}}]\label{P.product}
For $P$ and $Q$ finite posets,
$\Delta(P)\times\Delta(Q) \cong \Delta(P\times Q)$ as
partially ordered topological spaces, via the map
taking each element $(\,f_P,f_Q)\in\Delta(P)\times\Delta(Q)$
to the element $f\in\Delta(P\times Q)$ defined by
$f(x,y)=f_P(x)\wedge f_Q(y)$ $(x\in P,\ y\in Q,$
and the infimum ``$\wedge$\!'' taken in $[0,1]).$
\end{proposition}

\begin{proof}
That the function $f{:}\ P\times Q\to [0,1]$ given by
$f(x,y)=f_P(x)\wedge f_Q(y)$ belongs to $\Delta(P\times Q)$
is easily verified: for each $t\in\ho,$
$f_t$ will be $(\,f_P)_t\times(\,f_Q)_t,$ so writing
$(\,f_P)_t={\downarrow}\,h_t(\,f_P),$
$(\,f_Q)_t={\downarrow}\,h_t(\,f_Q),$ we see that
$f_t={\downarrow}\,h_t(\,f_P)\times {\downarrow}\,h_t(\,f_Q)
={\downarrow}\,(h_t(\,f_P),h_t(\,f_Q)),$
a principal ideal of $P\times Q.$

The resulting map $\Delta(P)\times\Delta(Q) \to \Delta(P\times Q)$ is
clearly isotone; i.e., if $(\,f_P,f_Q)\leq (g_P,g_Q)$
in $\Delta(P)\times\Delta(Q),$ then the elements of
$\Delta(P\times Q)$ constructed from these satisfy $f\leq g.$
Conversely, if $(\,f_P,f_Q)\not\leq (g_P,g_Q),$ then assuming
without loss of generality that $f_P\not\leq g_P,$ we can
choose $x\in P$ such that $f_P(x)>g_P(x).$
Taking a $y\in Q$ such that $f_Q(y)=1$ (i.e., any element
of the principal ideal $(\,f_Q)_1),$ we get
$f(x,y)=f_P(x)\wedge f_Q(y)=f_P(x)>g_P(x)\geq g_P(x)\wedge g_Q(y)=
g(x,y);$ so $f\not\leq g.$
Thus, the map described is an embedding of posets.

It remains to show that it is surjective.
Given $f\in\Delta(P\times Q),$ let us choose $(p,q)\in P\times Q$
such that $f(p,q)=1,$ and define $f_P(x)=f(x,q)$
and $f_Q(y)=f(p,y)$ for all $x\in P,$ $y\in Q.$
We want to show that $f_P$ and $f_Q$ belong to $\Delta(P)$
and $\Delta(Q)$ respectively, and that the element
$f'\in\Delta(P\times Q)$ defined by
\begin{equation}\begin{minipage}[c]{25pc}\label{d.f'}
$f'(x,y)\ =\ f_P(x)\wedge f_Q(y)$
\end{minipage}\end{equation}
is equal to $f.$
Thus, we need to show for every $t\in\ho$ that
the sets $(\,f_P)_t$ and $(\,f_Q)_t$ are principal ideals, and that
the sets $f'_t$ and $f_t$ coincide.

Given $t\in\ho,$ let
$f_t={\downarrow}\,(x_0,y_0)$ for $x_0\in P,$ $y_0\in Q.$
Since $f(p,q)=1,$ we have $(p,q)\in f_t = {\downarrow}\,(x_0,y_0),$
whence
\begin{equation}\begin{minipage}[c]{25pc}\label{d.pq}
$p\in {\downarrow}\,x_0$\quad and\quad $q\in {\downarrow}\,y_0.$
\end{minipage}\end{equation}
We now compute: $(\,f_P)_t=\{z\in P\mid f_P(z)\geq t\}
=\{z\in P\mid f(z,q)\geq t\}=\{z\in P\mid (z,q)\in f_t\}
=\{z\in P\mid (z,q)\in{\downarrow}\,(x_0,y_0)\}={\downarrow}\,x_0$
(the last step by the second relation of~\eqref{d.pq}).
Similarly $(\,f_Q)_t={\downarrow}\,y_0.$
So as required, $(\,f_P)_t$ and $(\,f_Q)_t$ are principal,
and by~\eqref{d.f'}, $f'_t={\downarrow}\,(x_0,y_0)=f_t.$
\end{proof}

\section{Stitching lattices together}\label{S.stitching}

\subsection{The construction}\label{S2.stitching}
We noted in \S\ref{S2.M3} that $\Delta(M_3)$ consists of
three triangles, each a distributive sublattice corresponding
to a maximal chain in $M_3,$ joined along a common edge.
Is this an instance of a general way that lattices can be attached
together?

In \cite[\S{}IV.2, esp.\ IV.2.3]{GG.Lat}, several ways of
attaching lattices to one another to get new
lattices are discussed, under the names
{\em gluing, pasting, S-gluing} and {\em multipasting}.
The next lemma shows that the way $\Delta(M_3)$ is put together
from sublattices can also be generalized.
We shall show in \S\ref{S2.not_gluing} that it
does not fall under any of those previously studied patterns.

\begin{lemma}\label{L.stitch}
Suppose $(P_i)_{i\in I}$ is a nonempty family of posets,
whose order-relations we shall write $\leq_i,$ such
that the $P_i$ are pairwise disjoint except for a common chain $C$
which has the same order structure in all the $P_i.$
Suppose moreover that for each $x\in P_i$ $(i\in I)$ there
is a least element $u_i(x)\in C$ that is $\geq_i x,$ and
a greatest element $d_i(x)\in C$ that is $\leq_i x.$
\textup{(}Note that if $C$ is {\em finite,} then this condition
holds if and only if each $P_i$ has a least and a greatest element,
and these lie in $C.$
The letters $u$ and $d$ are mnemonic for ``up'' and ``down''\textup{)}.

Then the set $P=\bigcup_I P_i$ may be partially ordered by taking
elements $x\in P_i,$ $y\in P_j$ to satisfy $x\leq y$
if and only if either $i=j$ and $x\leq_i y,$
or $i\neq j$ and the following equivalent conditions hold:
\textup{(i)}~$u_i(x)\leq d_j(y),$
\textup{(ii)}~$u_i(x)\leq_j y,$
\textup{(iii)}~$x\leq_i d_j(y),$
\textup{(iv)}~$x\leq_i z\leq_j y$ for some $z\in C.$

If the $P_i$ are $\!\vee\!$--semilattices, then the
resulting poset $P$ is also a $\!\vee\!$--semilattice;
namely, for $x\in P_i,$ $y\in P_j,$ we have
$x\vee y= x\vee_i y$ if $i=j,$ while if $i\neq j,$ and
we assume without loss of generality that $u_i(x)\geq u_j(y),$
then $x\vee y= x\vee_i u_j(y)\in P_i\subseteq P.$

If the $P_i$ are $\!\wedge\!$--semilattices, then
so is $P,$ by the dual construction.

Hence if the $P_i$ are lattices, $P$ is a lattice.
\end{lemma}

\begin{proof}[Sketch of proof]
The parenthetical observation on the case where $C$ is finite in the
first paragraph is immediate, as is the verification of
the equivalence of conditions~(i)-(iv) of the second paragraph.
The verification that the relation $\leq$ defined in that paragraph
is a partial ordering is routine; this includes the fact that it is
well-defined, namely, that if
$x$ and/or $y$ lies in $C,$ the condition for $x\leq y$ to hold
does not depend on which posets $P_i$ and/or $P_j$ those elements
are regarded as lying in.

Turning to the description of $x\vee y,$
it clearly gives (in each case) an
upper bound to $x$ and $y,$ so we must verify that this is
majorized by any upper bound $z$ of those elements, say lying in~$P_k.$

In the case $i=j,$ the desired result is clear if $k=i.$
If $k\neq i,$ then if $z$ majorizes $x$ and $y,$ it must
majorize $u_i(x)$ and $u_i(y),$ hence it majorizes
the larger of these, which is $\geq x\vee_i y,$ as required.

In the case where $i\neq j$ and $u_i(x)\geq u_j(y),$ we note
that $k$ must be distinct from at least one of $i$ and $j,$ hence
by definition of the ordering of $P,$
$z$ majorizes at least one of $u_i(x),$ $u_j(y);$ hence
it majorizes the smaller of them, $u_j(y).$
So in $P$ it majorizes both $x$ and $u_j(y),$
which both lie in $P_i;$ so by the preceding case
it majorizes $x\vee_i u_j(y),$ the asserted
join of $x$ and $y,$ as desired.

The dual assertion follows by symmetry.
Hence when the $P_i$ are lattices we get the final assertion.
\end{proof}

We remark that the assumption that $C$ is a chain was
not needed for our construction of the partial order
on $P,$ but only for the verification
that if the $P_i$ have meets and/or joins, so does $P;$
and for the parenthetical note on the case where $C$ is finite.

Following the theme of the terms ``gluing'' and
``pasting'', let us call the $P$ of
Lemma~\ref{L.stitch} the poset, semilattice,
or lattice obtained by {\em stitching} the $P_i$
together along $C.$

\subsection{When \texorpdfstring{$\Delta(P)$}{Delta(P)} can be obtained by stitching}\label{S2.stitch_Delta}
What instances of the stitching construction
can we hope will be respected by the construction $\Delta$?
For this to hold, not only $C,$ but also $\Delta(C)$ must be a chain;
hence $C$ must have $\leq 2$ elements.
Ignoring the trivial cases where it has zero or one element,
the next lemma shows that for $C$ having two elements
the desired result holds.
In particular, it gives
the description of $\Delta(M_3)$ which motivated these considerations.

\begin{lemma}\label{L.stitch_Delta}
Let $(P_i)_{i\in I}$ be a finite nonempty family of finite posets,
disjoint except for a common least element $0$
and a common greatest element $1,$ and let $P$ be the
poset obtained by stitching the $P_i$ together along $C=\{0,1\}.$

Then the poset $\Delta(P)$ can be obtained
by stitching together the $\Delta(P_i)$ along their
common $\!1\!$--simplex $\Delta(\{0,1\}).$
\end{lemma}

\begin{proof}
The case where the index-set $I$ has just one element is trivial, so
assume the contrary.
It is easy to see that for $P$ as described,
elements of distinct sets $P_i-C,$ $P_j-C$ will be
incomparable, hence every chain in $P$ must
be contained in one of the $P_i,$ hence will lie in more
than one of them if and only if it lies in $\{0,1\}.$
It follows from these observations, and the description of the operator
$\Delta$ in terms of chains
of elements, that $\Delta(P)$ can be identified
with the union of the $\Delta(P_i),$ and
that the $\Delta(P_i)$ are pairwise disjoint
except for the common $\!1\!$--simplex $\Delta(\{0,1\}).$

Note that for $x\in P_i,$ the principal ideal generated
by $x$ in $P$ is the same as the principal ideal generated
by $x$ within $P_i,$ except when $x=1,$
in which case it is all of $P.$
From this we can see that the embedding of $\Delta(P_i)$
in $\Delta(P)$ carries each $f\in\Delta(P_i)$ to the function
that agrees with $f$ on $P_i,$ and takes the value $f(1)$
everywhere else in $P.$
Likewise, the embeddings of $\Delta(\{0,1\})$ in the
spaces $\Delta(P_i)$ and $\Delta(P)$ take each
$f\in\Delta(\{0,1\})$ to the function which
has the value $f(1)$ at all elements other than $0.$

From this it follows that in $\Delta(P_i),$ each element $f$ is
majorized by
a smallest element $u_i(\,f)\in\Delta(\{0,1\}),$ namely, the function
whose value at $1$ is the greatest of the values of $f$
at points of $P_i$ other than $0.$
Likewise, $f$ majorizes a largest
element $d_i(\,f)\in\Delta(\{0,1\}),$ namely, the function
on $\{0,1\}$ whose value at $1$ is $f(1).$

It remains to verify that the order relation $f\leq g$ of $\Delta(P)$
is as described in the second paragraph of Lemma~\ref{L.stitch}.
We first note that for $f$ and $g$ in the same set $\Delta(P_i),$
our description of how elements of that
space extend to functions on $P$ shows
that they satisfy $f\leq g$ in $\Delta(P)$ if and only if they
satisfy the same inequality in $\Delta(P_i).$
From this case, it easily follows that for $f$ and $g$ lying in
distinct spaces $\Delta(P_i),$ $\Delta(P_j),$ the equivalent conditions
(i)-(iv) of Lemma~\ref{L.stitch} also imply that when $f$ and
$g$ are extended to elements of $\Delta(P),$ they satisfy
$f\leq g$ there.

Conversely, suppose $f\leq g$ in $\Delta(P),$
where $f\in\Delta(P_i)$ and $g\in\Delta(P_j)$ with $i\neq j.$
In particular, the values of $f$ at all points
of $P_i-\{0\}$ must be $\leq$ the values of $g$ at these points;
but the latter are $g(1);$ so $f$ is majorized by the function on $P$
which is $g(1)$ everywhere on $P-\{0\}.$
Since that function is $d_j(g),$ we get condition~(iii)
of Lemma~\ref{L.stitch}.
\end{proof}

\subsection{Stitching \texorpdfstring{$\!\notin \{\!$gluing etc.$\!\}$}{is not gluing, etc.}}\label{S2.not_gluing}
To verify that stitching of lattices does not fall under the
list of lattice constructions given in \cite[\S{}IV.2]{GG.Lat},
we recall that for all of those constructions,
the lattice constructed has the universal mapping property of
the colimit of the diagram formed from the given lattices
and sublattices.
(In the case of multipasting, this is made part of the
definition, \cite[Definition~6]{multipasting}.)
But the same is not true of stitching.
For instance, if we stitch together two three-element chains
$\{0<a<1\},$ $\{0<b<1\}$ along the common chain $\{0<1\},$
the result is a $\!4\!$--element lattice, isomorphic
to $\{0,1\}^2;$ but this does not have the universal
property referred to, since the inclusions of the two given lattices
in, say, the chain $\{0<a<b<1\}$ do not factor through $\{0,1\}^2.$

Note also that gluing and pasting preserve the class of
modular lattices \cite[Theorem~16]{pasting},
\cite[Theorem~303]{GG.Lat}; but stitching together the chains
$\{0<a<b<1\}$ and $\{0<c<1\}$ along $\{0,1\}$
gives the nonmodular lattice $N_5.$

\Needspace{5.5\baselineskip}
\section{Acknowledgements and notes}\label{S.akn}

I am indebted to George Gr\"{a}tzer, Patricia Hersh,
Fred Wehrung and Volkmar Welker
for invaluable pointers to what is known on some of these topics,
and to Walter Taylor for the same, and for the contributions noted
in the text, and, of course,
for his questions that inspired this investigation.
I am also indebted to
the referee for a long list of useful corrections and suggestions.

The arXiv number of this paper is 1602.00034\,.

After publication of this note, any updates, errata, related references
etc., if found, will be recorded at
\url{http://math.berkeley.edu/~gbergman/papers/}.

\end{document}